\documentclass{amsart}

\usepackage[usenames]{color}
\usepackage[dvips]{pstricks}
\usepackage[all]{xy}
\usepackage{amsmath,amsthm,amssymb,graphicx}
\usepackage{setspace}
\usepackage[normalem]{ulem}
 \theoremstyle{definition}
 
 \theoremstyle{remark}
 
 \numberwithin{equation}{subsection}
 \newcommand{\h}{\mathcal{H}}

\begin{document}

\newtheorem{thrm}{Theorem}
\newtheorem{lmm}[thrm]{Lemma}
\newtheorem{crllry}[thrm]{Corollary}
\newtheorem{cnjctr}[thrm]{Conjecture}
\newtheorem{dfntn}[thrm]{Definition}

\newcommand{\ignore}[1]{}
\newcommand{\cat}[1]{\mathcal{#1}}
\newcommand{\fld}[1]{\mathbb{#1}}
\newcommand{\leftexp}[2]{{\vphantom{#2}}^{#1}{#2}}
\newcommand{\standardfigsized}[3]{
  \begin{figure}
    \centerline{\includegraphics[#3]{#1.eps}}
\ignore{    \caption{#2 Reference #1.}\label{#1}     }
    \caption{#2}\label{#1}
  \end{figure}
}
\newcommand{\substandardfigsized}[3]{
  \begin{figure}
    \caption{INSERT PICTURE. #2 Reference #1.}\label{#1}
  \end{figure}
}
\newcommand{\standardfig}[2]{
  \standardfigsized{#1}{#2}{}
}\newcommand{\substandardfig}[2]{
  \substandardfigsized{#1}{#2}{}
}

\newcommand{\standardgraph}[1]{\includegraphics{#1.eps}}
\newcommand{\trivobject}{{\mathbf 1}}
\newcommand{\nnr}[2]{\langle #1,#2 \rangle}
\newcommand{\snnr}[2]{(#1,#2)}

\newcommand{\moon}[1]{{\color{blue}#1}}
\newcommand{\tobe}[1]{{\color{red}#1}}
\newcommand{\tobedel}[1]{{\color{red}\sout{#1}}}
\newcommand{\tobeadd}[1]{{\color{green}#1}}

\ignore{\newcommand{\fixme}[1]{}}

\title[Some non-braided fusion categories of rank 3]
 {Some non-braided fusion categories of rank 3}

\author{ Tobias J. Hagge and Seung-Moon Hong }

\address{Department of Mathematics,  Indiana University, Bloomington, Indana}

\email{thagge@indiana.edu,seuhong@indiana.edu}

\thanks{}

\thanks{}



\date{}

\dedicatory{}

\commby{}

%%% ----------------------------------------------------------------------

\begin{abstract}
We classify all fusion categories for a given set of fusion rules with three simple object types. If a conjecture of Ostrik is true, our classification completes the classification of fusion categories with three simple object types. To facilitate the discussion we describe a convenient, concrete and useful variation of graphical calculus for fusion categories, discuss pivotality and sphericity in this framework, and give a short and elementary re-proof of the fact that the quadruple dual functor is naturally isomorphic to the identity.

\end{abstract}

%%% ----------------------------------------------------------------------
\maketitle
%%% ----------------------------------------------------------------------

\section{Introduction}\label{sec_intro}
Let $k$ be an algebraically closed field. A {\em fusion category} $\cat{C}$ over $k$ is a $k$-linear semi-simple rigid monoidal category with finitely many (isomorphism classes of) simple objects, finite dimensional morphism spaces, and $End(\trivobject) \cong k$. See \cite{maclane} or \cite{kassel} for definitions, and \cite{etingof_nikshych_ostrik} for many of the known results about fusion categories.

The {\em rank} $r$ of $\cat{C}$ is the number of isomorphism classes of simple objects in $\cat{C}$. Let $\{x_i\}_{1 \le i \le r}$ be a set of simple object representatives. The {\em fusion rules} for $\cat{C}$ are a set of $r \times r$ $\fld{N}$-valued matrices $N=\{N_i\}_{1 \le i \le r}$, with $(N_i)_{j,k}$ denoted $N_{i j}^k$ or, when convenient, $N_{x_i x_j}^{x_k}$, such that $x_i \otimes x_j \cong \bigoplus_{1 \le k \le r} N_{i j}^k x_k$. In the sequel, assume $k=\fld{C}$.

Fusion categories appear in representation theory, operator algebras, conformal field theory, and in constructions of invariants of links, braids, and higher dimensional manifolds. There is currently no general classification of them. Classifications of fusion categories for various families of fusion rules have been given in work by Kerler (\cite{kerler}), Tambara and Yamagami (\cite{tambara_yamagami}), Kazhdan and Wenzl (\cite{kazhdan_wenzl}), and Wenzl and Tuba (\cite{wenzl_tuba}).

For a given set of fusion rules, there are only finitely many monoidal natural equivalence classes of fusion categories. This property is called Ocneanu rigidity (see \cite{etingof_nikshych_ostrik}). It is not known whether or not the number of fusion categories of a given rank is finite. If one assumes a modular structure, the possibilities up to rank four have been classified by Belinschi, Rowell, Stong and Wang in \cite{rowell_stong_wang}. Ostrik has classified fusion categories up to rank two in \cite{ostrik2}, and constructed a finite list of realizable fusion rules for braided categories up to rank three in \cite{qa0503564}, in which the number of categories for each set of fusion rules is known. The rank two classification relies in an essential way on the theory of modular tensor categories; Ostrik shows that the quantum double of a rank two category must be modular, and uses the theory of modular tensor categories to eliminate most of the possibilities. The classification of modular tensor categories is of independent interest; in many contexts one must assume modularity.

We consider the only set of rank three fusion rules which is known to be realizable as a fusion category but which has no braided realizations. Ostrik conjectured in \cite{qa0503564} that a classification for this rule set completes the classification of rank three fusion categories.

The axioms for fusion categories over $\fld{C}$ reduce to a system of polynomial equations over $\fld{C}$. In this context, Ocneanu rigidity, roughly translated, says that normalization of some of the variables in the equations gives a finite solution set. In this case, one can compute a Gr\"{o}bner basis for the system and obtain the solutions (see \cite{buchberger}). However, normalization becomes complicated when there are $i,j,k$ such that $N_{i j}^k > 1$. The fusion rules we consider are the smallest realizable set with this property.
\section{Main theorem and outline}\label{sec_outline}
\begin{thrm}[Main Theorem]\label{thm_main} Consider the set of fusion rules with three simple object types, $x$, $y$ and $\trivobject$. Let $\trivobject$ be the trivial object, and let $x \otimes x \cong x \oplus x \oplus y \oplus \trivobject$, $x \otimes y \cong y \otimes x \cong x$ and $y \otimes y \cong \trivobject$. Then the following hold:

\begin{enumerate}
\item \label{main_tensor}
Up to monoidal natural equivalence, there are four semisimple tensor categories with the above fusion rules. A set of associativity matrices for one of these categories is given in Appendix~\ref{sec_appendix}. Applying a nontrivial Galois automorphism to all of the coefficients gives a set of matrices for any one of the other three categories.
\item \label{main_rigidity} 
The categories in part \ref{main_tensor} are fusion categories.
\item \label{main_nobraid}
The categories in part \ref{main_tensor} do not admit braidings.
\item \label{main_spherical}
The categories in part \ref{main_tensor} are spherical.
\end{enumerate}
\end{thrm}

The structure of the remainder of the paper is as follows:

Section~\ref{sec_prelim} describes the notation and categorical preliminaries used in later parts of the paper. It constructs a canonical representative for each monoidal natural equivalence class of fusion categories. This construction is really just two well known constructions, skeletization and strictification, applied in sequence. These constructions, taken together, form a bridge between the category theoretic language in the statements of the theorems and the algebra appearing in the proofs. For some of the calculations in this paper the translation between the category theory and the algebra is already widely known, but there are some subtleties when discussing pivotal structure that justify the treatment. Section~\ref{sec_prelim} concludes by describing the algebraic equations corresponding to the axioms for a fusion category, using the language of strictified skeletons.

Section~\ref{sec_tensor} proves part \ref{main_tensor} of Theorem~\ref{thm_main}. The proof amounts to solving the variety of polynomial equations defined in the previous section, performing normalizations along the way in order to simplify calculations. The section ends by arguing that the nature of the normalizations guarantees that the solutions obtained really are monoidally inequivalent. Section~\ref{sec_rigidity} proves part \ref{main_rigidity} of Theorem~\ref{thm_main} by explicitly computing rigidity structures.

Part \ref{main_nobraid} of Theorem~\ref{thm_main} follows from Ostrik's classification of rank three braided fusion categories in \cite{qa0503564}. Section~\ref{sec_nobraid} gives a direct proof by showing that there are no solutions to the hexagon equations.

Section~\ref{sec_pivotalprelim} defines pivotal and spherical structures and discusses their properties. The focus is on the question of whether every fusion category is pivotal and spherical. A novel and elementary proof that the quadruple dual functor is naturally isomorphic to the identity functor is given. This proof makes use of the strictified skeleton construction developed in Section~\ref{sec_prelim}. The section concludes by describing what a pivotal category which does not admit a spherical structure would look like. In particular, it must have at least five simple objects.

Section~\ref{sec_pivotal} proves part \ref{main_spherical} of Theorem~\ref{thm_main} by computing explicit pivotal structures for the four categories given in Section~\ref{sec_tensor}, and invoking a lemma from Section~\ref{sec_pivotalprelim} for sphericity.

\section{Preliminaries and notational conventions}\label{sec_prelim}
This paper uses the ``composition of morphisms'' convention for functions as well as morphisms, and left to right matrix multiplication. For calculations of the fusion rules, our treatment is similar to \cite{tambara_yamagami}, but the notation differs superficially for typographic reasons. The notation captures algebraic data sufficient to classify a fusion category up to monoidal natural equivalence, and is reviewed later in this section.

We assume that the reader is familiar with the notions of a monoidal category, a monoidal functor and a monoidal equivalence; for precise definitions, see \cite{maclane}. Recall that a monoidal category is equipped with an associative bifunctor $\otimes$ and a distinguished object $\trivobject$. Reassociation of tensor factors in a monoidal category is described by a natural isomorphism of trifunctors $\alpha:(- \otimes -) \otimes - \to - \otimes (- \otimes -)$. Tensor products with $\trivobject$ have natural isomorphisms $\rho: - \otimes \trivobject \to -$ and $\lambda: \trivobject \otimes - \to -$. These isomorphisms are subject to a coherency condition, namely that for any pair of multifunctors there is at most one natural isomorphism between them which may be constructed from $\lambda$, $\rho$, $\alpha$ and their inverses, along with $Id$ and $\otimes$. This coherency condition is well known to be equivalent to the statement that the category satisfies the {\em pentagon} and {\em triangle} axioms (see \cite{maclane} for a proof).

The triangle equations are the equations $\rho_x \otimes y = \alpha_{x,1,y} \circ (x \otimes \lambda_y)$ for all ordered pairs $(x, y)$ of objects. Here, and in the sequel when the context is unambiguous, the name of an object is used as a shorthand for the identity morphism on that object. The pentagon equations, defined for all tuples of objects $(w, x, y, z)$, are as follows (see also Figure~\ref{fig:pentagon}):

$$
(\alpha_{x,y,z} \otimes w)\circ  \alpha_{x,y \otimes z,w} \circ( x
\otimes \alpha_{y,z,w})
 =  \alpha_{x
\otimes y,z,w}\circ \alpha_{x,y,z \otimes w},
$$

When studying fusion categories up to monoidal equivalence, one may choose categories within an equivalence class which have desirable attributes. Since the categorical properties considered in this paper (fusion rule structure, monoidality, pivotality, sphericity, presence of braidings) are all well known to be preserved under monoidal equivalence, the desirable attributes may be assumed without loss. In particular, one may construct, given an arbitrary fusion category, a canonical representative for that category's equivalence class in which one may replace instances of the words ``is isomorphic to'' with ``equals''.

\subsection{Skeletization}

The {\em skeleton} $\cat{C}^\cat{SKEL}$ of an arbitrary category $\cat{C}$ is any full subcategory of $\cat{C}$ containing exactly one object from each isomorphism class in $\cat{C}$. If $\cat{C}$ is semi-simple, every object in $\cat{C}$ is isomorphic to a direct sum of simple objects in $\cat{C}^\cat{SKEL}$. One may then assume without loss that the objects of $\cat{C}^\cat{SKEL}$ consist of simple object representatives and direct sums of such.

It is a well known fact that $\cat{C}^\cat{SKEL}$ may be given a monoidal structure such that $\cat{C}^\cat{SKEL}$ and $\cat{C}$ are monoidally equivalent. The proof is a straightforward but tedious extension of Maclane's proof of the natural equivalence an ordinary category and its skeleton (see section~IV.4 in \cite{maclane}). In that proof, one defines a family of isomorphisms $i_x$ from objects $x$ to their isomorphic representatives in $\cat{C}^\cat{SKEL}$ and uses it to construct a pair of functors $F$ and $G$ which give a natural equivalence. For the extension, the $i_x$ are used to define the tensor product functor on $\cat{C}^\cat{SKEL}$, as well as $\alpha$, $\lambda$, $\rho$ and the monoidal structures for $F$ and $G$. One then writes out all of the relevant commutative diagrams and removes any compositions $i_x \circ i_x^{-1}$. The result in each case is a commutative diagram in $\cat{C}$.

\subsection{Strictification}

Given a monoidal category $\cat{C}$, one may construct a strict monoidal category $\cat{C}^\cat{STR}$ equivalent to $\cat{C}$. In a strict monoidal category, $\alpha$, $\lambda$ and $\rho$ are the identity. It is common practice to assume that a monoidal category is strict without explicit reference to the construction. However, by using the construction explicitly we will be able to pick a canonical representative for an equivalence class of monoidal categories and provide a natural interpretation of the graphical calculus.

Strictification of a monoidal category is analogous to the construction of a tensor algebra; it gives an equivalent strict category $\cat{C}^\cat{STR}$ by replacing the tensor product with a strictly associative formal tensor product. The objects of $\cat{C}^\cat{STR}$ are finite sequences of objects in $\cat{C}$. Morphism spaces of the form
\[Mor((a_1,a_2, \ldots, a_{m-1},a_m), (b_1,b_2, \ldots, b_{n-1},b_n))\]
are given by
\[Mor(a_1 \otimes (a_2 \otimes \ldots (a_{m-1} \otimes a_m) \dots ),b_1 \otimes (b_2 \otimes \ldots (b_{n-1} \otimes b_n) \dots ))\]
in $\cat{C}$. The tensor product on objects is just concatenation of sequences, for morphisms it is the tensor product in $\cat{C}$ pre and post-composed with appropriate associativity morphisms. Monoidal equivalence of $\cat{C}$ with $\cat{C}^\cat{STR}$ is proven in section~XI.3 of \cite{maclane}.

It is not usually possible to make a fusion category strict and skeletal at the same time. However, the category $(\cat{C}^\cat{SKEL})^\cat{STR}$, while not a skeleton, is still unique up to strict natural equivalence. Also, it is a categorical realization of a graphical calculus, as will shortly become clear. The next subsection describes what strictified skeleta of fusion categories look like, up to strict equivalence.

\subsection{Strictified skeletal fusion categories}\label{sec_strictskel}

A {\em strictified skeletal fusion category} $\cat{C}$ is as follows: Let $N$ be a set of fusion rules for a set of objects $S$. Then the objects in $\cat{C}$ are multisets of finite sequences of elements of $S$. $\cat{C}$ has a tensor product $\otimes$, which is defined on objects by pairwise concatenation of sequences, distributed over elements of multisets. Direct sum of objects is given by multiset disjoint union.

 A {\em strand} is an object which is a sequence of length one. Strands correspond to simple object types. If $x$, $y$ and $z$ are strands, define $Mor(x \otimes y,z)$ to be a $k$ vector space isomorphic to $k^{N_{x y}^z}$. For brevity, $V_x^y$ will denote $Mor(x,y)$, and tensor products will be omitted when the context is clear.  A morphism is {\em $(n,m)$-stranded} if its source and target are sequences of length $n$ and $m$, respectively. A morphism is {\em $(n)$-stranded} if it is $(m,n-m)$-stranded for some $0 \le m \le n$.

 Semi-simplicity of $\cat{C}$ means that for all objects $w$, $x$, $y$ and $z$ there are vector space isomorphisms $\sum_{v \in S}V_{x y}^v \otimes V_{w v}^z \cong V_{w x y}^z \cong \sum_{v \in S} V_{w x}^v \otimes V_{v y}^z$. The first isomorphism is given by $f \otimes g \to (Id_w \otimes f) \circ g$, the inverse of the second by $h \otimes l \to (h \otimes Id_y) \circ l$. The composition of the two isomorphisms is denoted $\alpha_{w,x,y}^z$. Additionally, each morphism space $V_{x y}^z$ has an algebraically dual space $V_z^{x y}$, in the sense that there are bases $\{v_i\}_i \subset V_{x y}^z$ and $\{w_i\}_i \subset V_z^{x y}$ such that $w_i \circ v_j = \delta_{ij} Id_z$.

In a strictified skeleton, the trivial object $\trivobject$ is the zero length sequence. There is a strand which is isomorphic to the trivial object, but not equal. This strand shall be taken to be $\trivobject$ in the sequel. This choice makes the category non-strict, since $\lambda$ and $\rho$ are no longer the identity, but it is convenient for graphical calculus purposes.

(Right) rigidity in a strictified skeleton means that there is a set involution $*$ on the strands, and each strand $x$ has morphisms $b_x: \trivobject \to x \otimes x^*$ and $d_x: x^* \otimes x \to \trivobject$ such that $(Id_{x^*} \otimes b_{x}) \circ (d_x \otimes Id_{x^*}) = Id_{x^*}$ and $(b_x \otimes Id_{x}) \circ (Id_{x} \otimes d_x) = Id_{x}$. This implies that $N_{xy}^{z}=N_{z^*x}^{y^*}$. Left rigidity is similar, and in the sequel, the right rigidity morphism for $x^*$ will be defined to be the left rigidity morphism for $x$.

Define $*$, $b$, and $d$ on concatenations of strands such that $d_{x \otimes y}=(Id_{y^*} \otimes d_x \otimes Id_{y}) \circ d_y$ and extend to direct sums. Then there is a contravariant {\em (right) dual} functor $*$ which sends $f \in V_x^y$ to $f^* = (Id_{y^*} \otimes b_x) \circ (Id_{y^*} \otimes f \otimes Id_{x^*}) \circ (d_y \otimes Id_{x^*}) \in V_{y^*}^{x^*}$. The definition of a left dual functor is similar, and the two duals are inverse functors by rigidity.

Monoidality for a strictified skeletal fusion category implies that the $\alpha$ are identity morphisms. For this to be true, it is necessary and sufficient that the following equation holds for all objects $x, y, z, w$, and $u$. Each instance will be referred to as $P^u_{x, y, z, w}$ in the sequel.

$$
\bigoplus_{t} \alpha_{y,z,w}^{t} V^{u}_{x t} \circ
\bigoplus_{s}V^{s}_{y z} \alpha_{x,s,w}^{u} \circ \bigoplus_{t}
\alpha_{x,y,z}^{t} V^{u}_{tw}
 :
$$
$$
\bigoplus_{s,t}V^{s}_{z w} V^{t}_{y s} V^{u}_{xt}
 \rightarrow  \bigoplus_{s,t}V^{s}_{yz}
V^{t}_{sw} V^{u}_{xt} \rightarrow \bigoplus_{s,t}V^{s}_{yz}
V^{t}_{xs} V^{u}_{tw} \rightarrow \bigoplus_{s,t}V^{s}_{xy}
V^{t}_{sz} V^{u}_{tw}
$$

is equal to

$$
\bigoplus_{s}  V^{s}_{z w} \alpha_{x,y,s}^{u} \circ  \tau \circ
\bigoplus_{s}  V^{s}_{xy} \alpha_{s,z,w}^{u}
 :
$$
$$
\bigoplus_{s,t}V^{s}_{z w} V^{t}_{ys} V^{u}_{xt}
 \rightarrow  \bigoplus_{s,t}V^{s}_{zw}
V^{t}_{xy} V^{u}_{ts} \rightarrow \bigoplus_{s,t}V^{s}_{xy}
V^{t}_{zw} V^{u}_{st} \rightarrow \bigoplus_{s,t}V^{s}_{xy}
V^{t}_{sz} V^{u}_{tw}.
$$

Here $\otimes$ for vector spaces and morphisms are omitted, and $\tau$
is the isomorphism interchanging the first and the second factors of
vector space tensor products (see Figure~\ref{fig:pentagon}).

\begin{figure}
$$ \xymatrix{& {((x \otimes y) \otimes z) \otimes w} \ar[ld]_{\alpha_{x,y,z}w}
\ar[rd]^{\alpha_{xy,z,w}}\\ {(x \otimes (y \otimes z)) \otimes w}
\ar[d]_{\alpha_{x,yz,w}} & (a) & {(x \otimes y)\otimes(z \otimes w)}
\ar[d]^{\alpha_{x,y,zw}} \\ {x \otimes ((y \otimes z) \otimes w)}
\ar[rr]^{x \alpha_{y,z,w}} && {x \otimes (y \otimes (z \otimes
w))}}$$

$$ \xymatrix{\bigoplus_{s,t}V^{s}_{xy}
V^{t}_{sz} V^{u}_{tw}  && \bigoplus_{s,t}V^{s}_{xy} V^{t}_{z w}
V^{u}_{s t} \ar[ll]_{\bigoplus_{s}  V^{s}_{xy}
\alpha_{s,z,w}^{u}} \\
\bigoplus_{s,t}V^{s}_{yz} V^{t}_{xs} V^{u}_{tw}
\ar[u]^{\bigoplus_{t} \alpha_{x,y,z}^{t} V^{u}_{tw}} & (b) &
\bigoplus_{s,t}V^{s}_{zw}
V^{t}_{xy} V^{u}_{ts} \ar[u]_{\tau}\\
\bigoplus_{s,t}V^{s}_{yz} V^{t}_{sw} V^{u}_{xt}
\ar[u]^{\bigoplus_{s}V^{s}_{yz} \alpha_{x,s,w}^{u}} &&
\bigoplus_{s,t}V^{s}_{zw} V^{t}_{ys} V^{u}_{x
t}\ar[ll]_{\bigoplus_{t} \alpha_{y,z,w}^{t} V^{u}_{xt}}
\ar[u]_{\bigoplus_{s}  V^{s}_{zw} \alpha_{x,y,s}^{u}}}$$

 \caption{(a) Pentagon equality and (b) corresponding equality }
\label{fig:pentagon}
\end{figure}

\subsection{Remarks}

\begin{enumerate}
\item Every fusion category is monoidally naturally equivalent to a strictified skeleton. Also, two naturally equivalent strictified skeleta have an invertible equivalence functor that takes strands to strands. This implies that equivalences are given by permutations of strands along with changes of basis on the $(2,0)$ and $(2,1)$-stranded morphism spaces.

\item The functor $**$ fixes objects. The isomorphisms $J_{x,y}: x^{**} \otimes y^{**} \to (x \otimes y)^{**}$ associated with $**$ in the definition of a monoidal functor (see \cite{maclane}) may be taken to be trivial. There is an invertible scalar worth of freedom in the choice of each $b_x$, $d_x$ pair.

\item Semi-simplicity allows every morphism to be built up from $(3)$-stranded morphisms. Choosing bases for the $(3)$-stranded morphisms allows morphisms in $\cat{C}$ to be characterized as undirected trivalent graphs with labeled edges and vertices, subject to associativity relations given by the pentagon equations. The labels for the edges are isomorphism types of simple objects; the labels for the vertices are basis vectors for the corresponding morphism spaces. This gives a categorically precise interpretation of an arrowless graphical calculus for $\cat{C}$. 

\item If $\cat{C}$ is pivotal (see Section~\ref{sec_pivotalprelim} for the definition), a well known construction allows one to add a second copy of each object and get a strict pivotal category. This construction gives a graphical calculus with arrows on the strands.

\item Strictified skeleta give any categorical structure preserved under natural equivalence (and any functorial property preserved under natural isomorphism) a purely algebraic description.
\end{enumerate}
\section{Proof of Theorem~\ref{thm_main} part ~\ref{main_tensor}:possible tensor category structures}\label{sec_tensor}

In this section we classify, up to monoidal equivalence, all $\fld{C}$-linear semisimple tensor categories with fusion rules given in Theorem~\ref{thm_main}. This amounts to solving the matrix equations described in the previous section. The simplest equations (those involving $1 \times 1$ matrices) are solved first, and normalizations are performed as necessary in order to simplify the equations.

\subsection{Setting up the pentagon equations}
The fusion rules are given by
$x \otimes x \cong \trivobject \oplus y \oplus x \oplus x$, $x \otimes y \cong
y \otimes x \cong x$, and $y \otimes y \cong \trivobject$. The non-trivial
vector spaces are $V^{1}_{11}$, $V^{x}_{1x}$
,$V^{x}_{x1}$, $V^{y}_{1y}$, $V^{y}_{y1}$, $V^{x}_{xy}$,
$V^{x}_{yx}$, $V^{1}_{yy}$, $V^{1}_{xx}$, $V^{y}_{xx}$,and
$V^{x}_{xx}$, and they are all 1-dimensional except the last space
which is 2-dimensional.

Let's choose basis vectors in each space. If we fix any non-zero
vector $v^{1}_{11} \in V^{1}_{11}$, then there are unique vectors
$v^{x}_{1x} \in V^{x}_{1x}$, $v^{x}_{x1} \in
V^{x}_{x1}$,$v^{y}_{1y} \in V^{y}_{1y}$, and $v^{y}_{y1} \in
V^{y}_{y1}$ such that the triangle equality holds. For the other
spaces, choose any non-zero vectors in each space and denote them
by $v^{x}_{xy} \in V^{x}_{xy}$, $v^{x}_{yx} \in V^{x}_{yx}$,
$v^{1}_{yy} \in V^{1}_{yy}$, $v^{1}_{xx} \in V^{1}_{xx}$,
$v^{y}_{xx} \in V^{y}_{xx}$, $v_1$ and $v_2 \in V^{x}_{xx}$ where
the two vectors $v_1$ and $v_2$ are linearly independent.

There are 30 associativities. It is a well known fact that if at
least one of the bottom objects is $\trivobject$ then the
associativity is trivial. That is, with the above basis choices the matrix for $\alpha^{z}_{u,v,w}$ is trivial
if at least one of the $u,v$ and $w$ is $\trivobject$. Now we have
ten non-trivial 1-dimensional associativities,
$\alpha^{y}_{y,y,y}$,$\alpha^{x}_{x,y,y}$,$\alpha^{x}_{y,y,x}$,$\alpha^{1}_{x,y,x}$,
$\alpha^{y}_{x,y,x}$,$\alpha^{x}_{y,x,y}$,$\alpha^{1}_{x,x,y}$,$\alpha^{y}_{x,x,y}$,
$\alpha^{1}_{y,x,x}$, and $\alpha^{y}_{y,x,x}$, five non-trivial
2-dimensional ones,
$\alpha^{x}_{x,y,x}$,$\alpha^{x}_{x,x,y}$,$\alpha^{x}_{y,x,x}$,$\alpha^{1}_{x,x,x}$,
and $\alpha^{y}_{x,x,x}$, and one 6-dimensional one,
$\alpha^{x}_{x,x,x}$.

\subsection{Normalizations}

With the above basis choices we obtain a basis for each tensor product of vector spaces in a canonical way and can parameterize each
associativity and pentagon equation. However, at this point our basis elements have not been uniquely specified,
and we should expect to obtain solutions with free parameters. As
the calculation progresses it will be convenient to simplify the
pentagon equations by requiring certain coefficients of certain
associativity matrices to be $1$ or $0$. These normalizations should
be thought of as restrictions on the basis choices made above. Normalizations simplify the equations and have an additional
advantage: once the set of possible bases is sufficiently
restricted, Ocneanu rigidity \cite{etingof_nikshych_ostrik}
guarantees a finite set of possibilities for the associativity
matrices of fusion categories,
which can be found algorithmically by computing a Gr{\"o}bner basis.

\subsection{Associativity matrices}

The following are the 1-dimensional associativities:

$ \alpha^{y}_{y,y,y}:v^{1}_{yy}  v^{y}_{y1} \mapsto
a^{y}_{y,y,y}v^{1}_{yy} v^{y}_{1y} $

$ \alpha^{x}_{x,y,y}:v^{1}_{yy} v^{x}_{x1} \mapsto
a^{x}_{x,y,y}v^{x}_{xy} v^{x}_{xy} $

$ \alpha^{x}_{y,y,x}:v^{x}_{yx} v^{x}_{yx} \mapsto
a^{x}_{y,y,x}v^{1}_{yy} v^{x}_{1x} $

$ \alpha^{1}_{x,y,x}:v^{x}_{yx} v^{1}_{xx} \mapsto
a^{1}_{x,y,x}v^{x}_{xy} v^{1}_{xx} $

$ \alpha^{y}_{x,y,x}:v^{x}_{yx} v^{y}_{xx} \mapsto
a^{y}_{x,y,x}v^{x}_{xy} v^{y}_{xx} $

$ \alpha^{x}_{y,x,y}:v^{x}_{xy} v^{x}_{yx} \mapsto
a^{x}_{y,x,y}v^{x}_{yx} v^{x}_{xy} $

$ \alpha^{1}_{x,x,y}:v^{x}_{xy} v^{1}_{xx} \mapsto
a^{1}_{x,x,y}v^{y}_{xx} v^{1}_{yy} $

$ \alpha^{y}_{x,x,y}:v^{x}_{xy} v^{y}_{xx} \mapsto
a^{y}_{x,x,y}v^{1}_{xx} v^{y}_{1y} $

$ \alpha^{1}_{y,x,x}:v^{y}_{xx} v^{1}_{yy} \mapsto
a^{1}_{y,x,x}v^{x}_{yx} v^{1}_{xx} $

$ \alpha^{y}_{y,x,x}:v^{1}_{xx} v^{y}_{y1} \mapsto
a^{y}_{y,x,x}v^{x}_{yx} v^{y}_{xx} $

where associativity coefficients are all non-zero.

For 2-dimensional and 6-dimensional associativities we need to fix
the ordering of basis elements in each Hom vector space. The orderings are as follows:

$\{ v^{x}_{yx} v_1, v^{x}_{yx} v_2 \}$ for $V^{x}_{x(yx)}$,
$\{
v^{x}_{xy} v_1, v^{x}_{xy} v_2 \}$ for $V^{x}_{(xy)x}$,

$\{
v^{x}_{xy} v_1, v^{x}_{xy} v_2 \}$ for $V^{x}_{x(xy)}$,
$\{
 v_1 v^{x}_{xy}, v_2 v^{x}_{xy}  \}$ for $V^{x}_{(xx)y}$,

$\{
v_1 v^{x}_{yx} , v_2 v^{x}_{yx}  \}$ for $V^{x}_{y(xx)}$,
$\{
v^{x}_{yx} v_1, v^{x}_{yx} v_2 \}$ for $V^{x}_{(yx)x}$,

$\{
 v_1 v^{1}_{xx}, v_2 v^{1}_{xx}  \}$ for $V^{1}_{x(xx)}$,
 $\{ v_1 v^{1}_{xx}, v_2 v^{1}_{xx}  \}$ for $V^{1}_{(xx)x}$,

 $\{ v_1 v^{y}_{xx}, v_2 v^{y}_{xx}  \}$ for $V^{y}_{x(xx)}$,
 $\{ v_1 v^{y}_{xx}, v_2 v^{y}_{xx}  \}$ for $V^{y}_{(xx)x}$,

 $\{ v^{1}_{xx} v^{x}_{x1}, v^{y}_{xx} v^{x}_{xy}, v_1 v_1, v_1 v_2, v_2 v_1, v_2 v_2 \}$ for
 $V^{x}_{x(xx)}$,

and
$\{ v^{1}_{xx} v^{x}_{1x}, v^{y}_{xx} v^{x}_{yx}, v_1 v_1, v_1 v_2,
v_2 v_1, v_2 v_2 \}$ for
 $V^{x}_{(xx)x}$.

With these ordered bases, each associativity has a matrix form
 (recall that we are using the right multiplication convention). That is,
 $\alpha^{x}_{x,y,x}$ is given by the invertible $2 \times 2$ matrix
 $a^{x}_{x,y,x}$, and $\alpha^{x}_{x,x,y}$ is given by the invertible $2 \times 2$ matrix
 $a^{x}_{x,x,y}$,etc., and finally $\alpha^{x}_{x,x,x}$ is given by the invertible $6 \times 6$ matrix
 $a^{x}_{x,x,x}$.

\subsection{Pentagon equations with $1 \times 1$ matrices}

Considering only nontrivial associativities, there are 17
1-dimensional pentagon equations, 14 2-dimensional pentagon
equations, 6 6-dimensional ones, and 1 16-dimensional one. Without
redundancy, the following are the 1-dimensional equations:

$P^{x}_{x,y,y,y}$ : $a^{y}_{y,y,y}$ $a^{x}_{x,y,y}=$
$a^{x}_{x,y,y}$.

$P^{1}_{x,x,y,y}$ : $a^{x}_{x,y,y}$$a^{1}_{x,x,y}$ $a^{y}_{x,x,y}=1$
,

$P^{1}_{x,y,x,y}$ : $a^{x}_{y,x,y}$ $a^{y}_{x,y,x}=$ $a^{1}_{x,y,x}$
,

$P^{y}_{x,y,x,y}$ : :$a^{x}_{y,x,y}$ $a^{1}_{x,y,x}=$
$a^{y}_{x,y,x}$ ,

$P^{1}_{x,y,y,x}$ : $a^{x}_{y,y,x}$ $a^{x}_{x,y,y}=$
$(a^{1}_{x,y,x})^2$ ,

$P^{y}_{x,y,y,x}$ :
 $a^{x}_{y,y,x}$ $a^{x}_{x,y,y}=$ $(a^{y}_{x,y,x})^2$ ,

$P^{x}_{y,x,y,y}$ :$(a^{x}_{y,x,y})^2=1$ ,

$P^{1}_{y,y,x,x}$ : $a^{y}_{y,x,x}$ $a^{1}_{y,x,x}$
$a^{x}_{y,y,x}=1$ ,

$P^{1}_{y,x,x,y}$ : $a^{y}_{x,x,y}$ $a^{y}_{y,x,x}=$ $a^{1}_{y,x,x}$
$a^{1}_{x,x,y}$

If we normalize the basis we may assume $a^{x}_{y,y,x},
a^{1}_{x,y,x}$ and $a^{1}_{x,x,y}$ to be 1 (for normalization see
\cite{tambara_yamagami} or \cite{kitaev}), and we can solve the
above 1-dimensional equations. Here is the solution:

$a^{y}_{y,y,y}=$ $a^{x}_{x,y,y}=$ $a^{y}_{x,x,y}=1$,
$a^{x}_{y,x,y}=$ $a^{y}_{x,y,x}=\pm 1$, $a^{1}_{y,x,x}=$
$a^{y}_{y,x,x}= \pm 1$.

Let's say $g:=a^{x}_{y,x,y}=a^{y}_{x,y,x}$ and
$h:=a^{1}_{y,x,x}=a^{y}_{y,x,x}$ in the sequel. Also let
$A:=a^x_{x,y,x}$, $B:=a^x_{x,x,y}$, $D:=a^1_{x,x,x}$,
$E:=a^y_{x,x,x}$, $F:=a^x_{y,x,x}$ and $\Phi:=a^x_{x,x,x}$ for
brevity.

\subsection{Pentagon equations with $2 \times 2$ or $6 \times 6$ matrices}

Now, the following are the 2-dimensional pentagon equations using
the above 1-dimensional solutions:

$P^{x}_{y,y,x,x}$ : $F^2=Id_2$

$P^{x}_{y,x,y,x}$ : $g A F=FA$

$P^{x}_{x,y,y,x}$ : $A^2=Id_2$

$P^{x}_{y,x,x,y}$ : $g B F=FB$

$P^{x}_{x,y,x,y}$ : $g B A=AB$

$P^{x}_{x,x,y,y}$
 : $B^2=Id_2$

$P^{1}_{y,x,x,x}$ :  $E F=D$

$P^{y}_{y,x,x,x}$ : $D F=E$

$P^{1}_{x,y,x,x}$
 : $F D A=D$

$P^{y}_{x,y,x,x}$
 : $F E A= g E$

$P^{1}_{x,x,y,x}$
 : $A D B= D$

$P^{y}_{x,x,y,x}$
 : $A EB= g E$

$P^{1}_{x,x,x,y}$ : $B E =D$

$P^{y}_{x,x,x,y}$ : $B D = E$

It should be noted that for this particular category the large
number of one dimensional morphism spaces gives us $q$-commutativity
relations and matrices with $\pm 1$ eigenvalues, which are of great
help when simplifying the pentagon equations by hand.

To analyze 2-dimensional and 6-dimensional pentagon equations, at first let's look at
the isomorphism $\tau $ interchanging the first and the second
factors of tensor products. This change of basis is necessary for 6-dimensional pentagon equations
because the image basis of the matrix for $\alpha_{x,y,zw}^u$ and
the domain basis of the matrix for $\alpha_{xy,z,w}^{u}$ may not be
the same. For $P^{x}_{x,y,x,x}$, $\tau $ is an isomorphism from the
space $V^1_{xx} V^x_{xy} V^x_{x1} \oplus V^y_{xx} V^x_{xy}
V^x_{xy}\oplus V^x_{xx}V^x_{xy}V^x_{xx}$ to $V^x_{xy} V^1_{xx}
V^x_{x1} \oplus V^x_{xy} V^y_{xx} V^x_{xy}\oplus
V^x_{xy}V^x_{xx}V^x_{xx}$, both of which correspond to $Hom
((x\otimes y)\otimes(x\otimes x),x)$. With the canonically ordered
basis $\{v^1_{xx}v^x_{xy}v^x_{x1}, v^y_{xx}v^x_{xy}v^x_{xy},
v_iv^x_{xy}v_j \}$ and $\{v^x_{xy}v^1_{xx}v^x_{x1}$,
$v^x_{xy}v^y_{xx}v^x_{xy}$, $v^x_{xy}v_iv_j \}$, respectively,
$\tau$ turns out to be $I_6$. For $P^{x}_{y,x,x,x}$,
$P^{x}_{x,x,y,x}$ and $P^{x}_{x,x,x,y}$, $\tau$ is also $I_6$. But
for $P^{1}_{x,x,x,x}$, it is $\tau_1$, and for $P^{y}_{x,x,x,x}$, it
is $\tau_2$, defined as follows:

 $\tau_1:= \left[\begin{array}{cccccc} 1 & 0 & 0 & 0 & 0 & 0
\\ 0&1&0&0&0&0\\0&0&1&0&0&0\\
0&0&0&0&1&0\\0&0&0&1&0&0\\0&0&0&0&0&1  \end{array} \right]$,
 $\tau_2:=\left[\begin{array}{cccccc} 0&1&0&0&0&0\\ 1&0&0&0&0&0\\0&0&1&0&0&0\\
0&0&0&0&1&0\\0&0&0&1&0&0\\0&0&0&0&0&1  \end{array} \right]$.

Here are the 6 6-dimensional pentagon equations:

$P^{x}_{y,x,x,x}$ ;  $ \Phi (I_2 \oplus I_2 \otimes F) \left( \bigl[
\begin{smallmatrix}0&h\\h&0 \end{smallmatrix} \bigr] \oplus F\otimes I_2
\right)=\left(\bigl[
\begin{smallmatrix}1&0\\0&g\end{smallmatrix}
\bigr] \oplus I_2 \otimes F \right) \Phi $

$P^{x}_{x,y,x,x}$ ; $\left( \bigl[
\begin{smallmatrix}0&h\\h&0\end{smallmatrix}
\bigr] \oplus F \otimes I_2\right) \Phi \left(\bigl[
\begin{smallmatrix}1&0\\0&g\end{smallmatrix}
\bigr] \oplus A \otimes I_2\right)= (I_2 \oplus I_2 \otimes A)\Phi
$

$P^{x}_{x,x,y,x}$ ; $\left(\bigl[
\begin{smallmatrix}1&0\\0&g\end{smallmatrix}
\bigr] \oplus A\otimes I_2\right) \Phi \left(\bigl[
\begin{smallmatrix}0&1\\1&0\end{smallmatrix}
\bigr] \oplus B\otimes I_2\right)= \Phi(I_2 \oplus I_2 \otimes A)  $

$P^{x}_{x,x,x,y}$ ; $\left(\bigl[
\begin{smallmatrix}0&1\\1&0\end{smallmatrix}
\bigr] \oplus B \otimes I_2 \right) ( I_2 \oplus I_2 \otimes B) \Phi
= \Phi\left( \bigl[
\begin{smallmatrix}1&0\\0&g\end{smallmatrix}
\bigr] \oplus I_2 \otimes B\right)  $

$ P^{1}_{x,x,x,x}$ ; $ (I_2 \oplus I_2 \otimes D)
 \Phi =( I_2 \oplus I_2 \otimes D) \tau_1
\left(\bigl[
\begin{smallmatrix}1&0\\0&h\end{smallmatrix}
\bigr] \oplus I_2 \otimes D\right)  $

$P^{y}_{x,x,x,x}$ ; $\Phi\left(\bigl[
\begin{smallmatrix}1&0\\0&g\end{smallmatrix}
\bigr] \oplus I_2 \otimes E\right)
 \Phi =( I_2 \oplus I_2 \otimes E) \tau_2
\left(\bigl[
\begin{smallmatrix}1&0\\0&h\end{smallmatrix}
\bigr] \oplus I_2 \otimes E\right) $.

If we normalize the basis $\{v_1,v_2 \}$ of $V^x_{xx}$, we may
assume $A$ is of the form $\bigl[
\begin{smallmatrix}1&0\\0&-1\end{smallmatrix}
\bigr]$, and then get $g=-1$ and $h=1$ using the above equations.
Following is the computation for this:

At first we may assume that matrix $A$ is of the Jordan canonical
form, then $A=$  $\pm I_2$ or $ \bigl[
\begin{smallmatrix} 1&0 \\0&-1 \end{smallmatrix}
\bigr]$ from $P^{x}_{x,y,y,x}$. We eliminate the possibility $A=$
$\pm I_2$ from $P^{x}_{y,x,y,x}$, $P^{1}_{x,y,x,x}$ and
$P^{x}_{y,x,x,x}$ which imply respectively that $g=1$, $F=\pm I_2$ and then $det( \Phi
)=0$, since the first two columns of $\Phi$ are scalar multiples of each other. So we conclude $A=$ $ \bigl[
\begin{smallmatrix} 1&0 \\0&-1 \end{smallmatrix}
\bigr]$. Now we eliminate the possibility $g=1$ using
$P^{x}_{y,x,y,x}$, $P^{x}_{y,y,x,x}$ and $P^{x}_{y,x,x,x}$, which
imply $F$ is a diagonal matrix, with entries $\pm 1$ and then $ det
( \Phi ) =0$, respectively. For the case $A=$ $\bigl[
\begin{smallmatrix}1&0\\0&-1\end{smallmatrix}
\bigr]$ and $g=-1$, $F$ is of the form $\bigl[
\begin{smallmatrix}0&f\\1/f&0\end{smallmatrix}
\bigr]$ from $P^{x}_{y,x,y,x}$ and $P^{x}_{y,y,x,x}$, and $B$ is of
the form $\bigl[
\begin{smallmatrix}0&b\\1/b&0\end{smallmatrix}
\bigr]$ from $P^{x}_{x,y,x,y}$ and $P^{x}_{x,x,y,y}$. If $h=-1$, the
first column of $\Phi$ has to be zero by comparing the first and the
second columns of $P^{x}_{y,x,x,x}$, $P^{x}_{x,y,x,x}$,
$P^{x}_{x,x,y,x}$, and $P^{x}_{x,x,x,y}$.

At this point we have fixed all 1-dimensional associativity matrices.

From the above equations, we get $\bigl[
\begin{smallmatrix}0&f\\1/f&0\end{smallmatrix}
\bigr]$ for $F$ and $\bigl[
\begin{smallmatrix}0&b\\1/b&0\end{smallmatrix}
\bigr]$ for $B$ with the relation $f^2 + b^2 =0$ from
$P^{x}_{y,x,x,y}$. We note that the diagonalization of $A$ defines each basis element $v_1$ and $v_2$ only up to choice of a nonzero scalar. By using up one of these degrees of freedom, we may
assume $f=1$. Then from the above 6-dimensional equations, we get
the following:

$D=  d \left[
\begin{array}{cc}1&b\\1&-b\end{array} \right]$

$E=  d \left[
\begin{array}{cc}b&1\\-b&1\end{array} \right]$

$\Phi=   \left[
\begin{array}{cccccc} \phi & \phi & -wb & w&w&-wb\\
\phi & -\phi & -wb & w&-w&wb \\x&x&-yb&z&y&-zb \\ x&x&-zb&y&z&-yb\\
x&-x&-yb&z&-y&zb\\ -x&x&zb&-y&z&-yb
\end{array} \right]$

\subsection{The pentagon equation with $16 \times 16$ matrices}

Now we analyze the 16-dimensional pentagon equation
$P^{x}_{x,x,x,x}$. It is convenient to express each Hom vector space
in two different ways and put basis permutation matrices into the
pentagon equation. The following are two expressions with ordered
direct sum.

$Hom(x(x(xx)),x):$

$\:\:\: V^x_{xx} V^1_{xx} V^x_{x1} \oplus V^x_{xx} V^y_{xx} V^x_{xy}
\oplus V^1_{xx} V^x_{x1} V^x_{xx} \oplus V^y_{xx} V^x_{xy} V^x_{xx}
\oplus V^x_{xx}V^x_{xx} V^x_{xx}$, and

$\:\:\:V^1_{xx} V^x_{x1} V^x_{xx} \oplus V^y_{xx} V^x_{xy} V^x_{xx}
\oplus V^x_{xx} V^1_{xx} V^x_{x1} \oplus V^x_{xx} V^y_{xx} V^x_{xy}
\oplus V^x_{xx}V^x_{xx} V^x_{xx}$

$Hom(x((xx)x)),x)$ :

$\:\:\:V^x_{xx} V^1_{xx} V^x_{x1} \oplus V^x_{xx} V^y_{xx} V^x_{xy}
\oplus V^1_{xx} V^x_{1x} V^x_{xx} \oplus V^y_{xx} V^x_{yx} V^x_{xx}
\oplus V^x_{xx}V^x_{xx} V^x_{xx}$, and

 $\:\:\: V^1_{xx} V^x_{1x} V^x_{xx}
\oplus V^y_{xx} V^x_{yx} V^x_{xx} \oplus V^x_{xx} V^1_{xx} V^x_{x1}
\oplus V^x_{xx} V^y_{xx} V^x_{xy} \oplus V^x_{xx}V^x_{xx} V^x_{xx}$

$Hom((x(xx))x),x)$ :

$\:\:\:V^1_{xx} V^x_{x1} V^x_{xx} \oplus V^y_{xx} V^x_{xy} V^x_{xx}
\oplus V^x_{xx} V^1_{xx} V^x_{1x} \oplus V^x_{xx} V^y_{xx} V^x_{yx}
\oplus V^x_{xx}V^x_{xx} V^x_{xx}$, and

$ \:\:\:V^x_{xx} V^1_{xx} V^x_{1x} \oplus V^x_{xx} V^y_{xx} V^x_{yx}
\oplus V^1_{xx} V^x_{x1} V^x_{xx} \oplus V^y_{xx} V^x_{xy} V^x_{xx}
\oplus V^x_{xx}V^x_{xx} V^x_{xx}$

$Hom((((xx)x)x),x)$ :

$\:\:\:V^x_{xx} V^1_{xx} V^x_{1x} \oplus V^x_{xx} V^y_{xx} V^x_{yx}
\oplus V^1_{xx} V^x_{1x} V^x_{xx} \oplus V^y_{xx} V^x_{yx} V^x_{xx}
\oplus V^x_{xx}V^x_{xx} V^x_{xx}$, and

$\:\:\:V^1_{xx} V^x_{1x} V^x_{xx} \oplus V^y_{xx} V^x_{yx} V^x_{xx}
\oplus V^x_{xx} V^1_{xx} V^x_{1x} \oplus V^x_{xx} V^y_{xx} V^x_{yx}
\oplus V^x_{xx}V^x_{xx} V^x_{xx}$

$Hom((xx)(xx),x)$ :

$\:\:\:V^1_{xx} V^x_{xx} V^x_{x1} \oplus V^y_{xx} V^x_{xx} V^x_{xy}
\oplus V^x_{xx} V^1_{xx} V^x_{1x} \oplus V^x_{xx} V^y_{xx} V^x_{yx}
\oplus V^x_{xx}V^x_{xx} V^x_{xx}$, and

$\:\:\:  V^1_{xx}V^x_{xx} V^x_{1x} \oplus V^y_{xx} V^x_{xx} V^x_{yx}
\oplus V^x_{xx} V^1_{xx} V^x_{x1} \oplus V^x_{xx} V^y_{xx} V^x_{xy}
\oplus V^x_{xx}V^x_{xx} V^x_{xx}$

 where each direct summand space has canonical ordered basis. For example
 $V^x_{xx}V^x_{xx} V^x_{xx}$ has basis $ \{v_i v_j v_k \}$ where $(i,j,k)$ range from 1 to 2
 in the order $(1,1,1),(1,1,2),(1,2,1),$ etc., and $V^x_{xx} V^1_{xx} V^x_{1x}$
 has $\{v_1 v^1_{xx} v^x_{1x},v_2 v^1_{xx} v^x_{1x} \}$.

Let

$\tau_3:= \left( \left[
\begin{array}{cc}0&1\\1&0\end{array} \right] \otimes I_4 \right) \oplus I_8$
and

$\tau_4:= \left( \left[
\begin{array}{cc}0&1\\1&0\end{array} \right] \otimes I_4 \right) \oplus
\left( \left[
\begin{array}{cccc}1&0&0&0\\0&0&1&0\\0&1&0&0\\0&0&0&1 \end{array}
\right] \otimes I_2 \right)$.

Then the pentagon equation $P^{x}_{x,x,x,x}$ is of the form:

$$ ( D \oplus E \oplus \Phi \otimes I_2) \tau_3 (I_2 \oplus A \oplus
\tilde{\Phi}) \tau_3 ( D \oplus E \oplus \Phi \otimes I_2) \tau_3
$$
$$  =\tau_3 (I_2 \oplus B \oplus \tilde{\Phi}) \tau_4 (I_2 \oplus F
\oplus \tilde{\Phi} )$$

where

$\tilde{\Phi}=   \left[
\begin{array}{cccccccccccc}
\phi &0& \phi&0 & -wb & w&w&-wb&0&0&0&0
\\ 0& \phi &0& \phi &0&0&0&0& -wb & w&w&-wb
\\ \phi &0& -\phi &0& -wb & w&-w&wb &0&0&0&0
\\0&\phi &0& -\phi &0&0&0&0& -wb & w&-w&wb
\\x&0&x&0&-yb&z&y&-zb&0&0&0&0
\\ x&0&x&0&-zb&y&z&-yb&0&0&0&0
\\x&0&-x&0&-yb&z&-y&zb&0&0&0&0
\\ -x&0&x&0&zb&-y&z&-yb&0&0&0&0
\\0&x&0&x&0&0&0&0&-yb&z&y&-zb
\\ 0&x&0&x&0&0&0&0&-zb&y&z&-yb
\\0&x&0&-x&0&0&0&0&-yb&z&-y&zb
\\ 0&-x&0&x&0&0&0&0&zb&-y&z&-yb
\end{array} \right]$.

\subsection{Solutions}

We may assume $x=1$ once we normalize basis vector $v^x_{xy}$.
Then from the equations, we get four explicit solution sets for
the parameters $b,\phi,d,w,y$, and $z$. We list one solution here. All of its values lie in the field $\fld{Q}(\sqrt{3},i)$; the other solutions are obtained by applying Galois automorphisms. The full set of associativity matrices for this solution is given in Appendix~\ref{sec_appendix}.

\[b = i,  \phi = \frac{-1+\sqrt{3}}{2}, d = \frac{1}{\sqrt{2}} e^{7 \pi i
 /12},   w = \frac{1- \sqrt{3}}{4} e^{2\pi i /3},\]
\[  y = \frac{1}{2} (e^{-\pi i /3} + i), z = \frac{1}{2} e^{5\pi i /6}.
 \]

\subsection{Inequivalence of the solutions}

To see that these solutions are monoidally inequivalent, recall from the previous section that for strictified skeletons a natural equivalence between two solutions to the pentagon equations is limited to change of basis on the $(2,0)$ and $(2,1)$-stranded morphism spaces, along with permutation of the strands. In our case permutation of strands does not preserve the fusion rules. Therefore, we must show that it is not possible to replicate the effect of a nontrivial Galois automorphism by change of basis choices for the $(2,0)$ and $(2,1)$-stranded morphism spaces.

The Galois automorphism that fixes $\sqrt{3}$ and sends $i$ to $-i$ changes the eigenvalues of the matrix $a^1_{x,x,x}$. However, $a^1_{x,x,x}$ is determined by the basis choices $v_1$ and $v_2$, and its rows and columns are indexed by $v_1$ and $v_2$. Thus changes to $v_1$ and $v_2$ conjugate $a^1_{x,x,x}$ by a change of basis matrix, which doesn't affect its eigenvalues. Therefore this automorphism does not correspond to a change of basis.

The other two Galois automorphisms send $\sqrt{3}$ to $-\sqrt{3}$ and thus change the value of the $(1,1)$ entry of $a^x_{x,x,x}$. But this entry is invariant under change of basis. Therefore no Galois automorphism corresponds to a change in basis, and the four solutions given above are mutually monoidally inequivalent.

\section{Proof of Theorem~\ref{thm_main} part \ref{main_rigidity}: rigidity structures}\label{sec_rigidity}

This section explicitly computes rigidity structures for the categories given in the previous section. Rigidity implies that these categories are fusion categories.

Given $v^1_{xx}\in V^1_{xx}$, choose a vector $ v^{x x}_1\in V^{x
x}_1$ such that $ v^{xx}_1 \circ v^1_{x,x} = id_1$(see Figure
\ref{fig:basis}). Now we define right death and birth, $d_x :=
v^1_{xx} :x \otimes x \rightarrow \trivobject,$ $b_x:=
\frac{1}{\phi} v^{xx}_1: \trivobject \rightarrow x \otimes x$ (see
Figure \ref{fig:b and d}).

With these definitions, right rigidity is an easy consequence by
direct computation. The following is a graphical version of it:

\psset{unit=3mm}
\begin{pspicture}[.3](0,0)(3,5)

\psline{<-}(0,5)(0,2.5) \psbezier(0,2.5)(0,0)(1.5,0)(1.5,2.5)
\psbezier(1.5,2.5)(1.5,5)(3,5)(3,2.5) \psline(3,2.5)(3,0)
\end{pspicture}
$= \frac{1}{\phi}$
\begin{pspicture}[.3](0,0)(5,7)
\psline(2,7)(2,6)
\psline(2,6)(1,5)\psline[linestyle=dotted](2,6)(3,5)
\psline(1,5)(0,4)\psline(3,5)(2,4)\psline(3,5)(4,4)
\psline(0,4)(0,3)\psline(2,4)(2,3)\psline(4,4)(4,3)
\psline(0,3)(1,2)\psline(2,3)(1,2)\psline(4,3)(3,2)
\psline(3,2)(2,1) \psline(2,1)(2,0)
\end{pspicture}
$=$
\begin{pspicture}[.3](0,0)(5,7)
\psline(2,7)(2,6) \psline(2,6)(3,5)
\psline(1,5)(0,4)\psline(1,5)(2,4)\psline(3,5)(4,4)
\psline(0,4)(0,3)\psline(2,4)(2,3)\psline(4,4)(4,3)
\psline(0,3)(1,2)\psline(2,3)(1,2)\psline(4,3)(3,2)
\psline(3,2)(2,1) \psline(2,1)(2,0)
\end{pspicture}
$=$
\begin{pspicture}[.3](0,0)(3,5)
\psline(1,5)(1,0)
\end{pspicture} $= id_x$,

\psset{unit=3mm}
\begin{pspicture}[.3](0,0)(3,5)
\psline{<-}(0,0)(0,2.5) \psbezier(0,2.5)(0,5)(1.2,5)(1.5,2.5)
\psbezier(1.5,2.5)(1.5,0)(3,0)(3,2.5) \psline(3,2.5)(3,5)
\end{pspicture}
$= \frac{1}{\phi}$
\begin{pspicture}[.3](0,0)(5,7)
\psline(2,7)(2,6)
\psline[linestyle=dotted](2,6)(1,5)\psline(2,6)(3,5)
\psline(1,5)(0,4)\psline(1,5)(2,4)\psline(3,5)(4,4)
\psline(0,4)(0,3)\psline(2,4)(2,3)\psline(4,4)(4,3)
\psline(0,3)(1,2)\psline(2,3)(3,2)\psline(4,3)(3,2)
\psline(1,2)(2,1) \psline(2,1)(2,0)
\end{pspicture}
$=$
\begin{pspicture}[.3](0,0)(5,7)
\psline(2,7)(2,6) \psline(2,6)(1,5)
\psline(1,5)(0,4)\psline(3,5)(2,4)\psline(3,5)(4,4)
\psline(0,4)(0,3)\psline(2,4)(2,3)\psline(4,4)(4,3)
\psline(0,3)(1,2)\psline(2,3)(3,2)\psline(4,3)(3,2)
\psline(1,2)(2,1) \psline(2,1)(2,0)
\end{pspicture}
$=$
\begin{pspicture}[.3](0,0)(3,5)
\psline(1,5)(1,0)
\end{pspicture}$= id_x$

where the first and the third equalities are from the definitions
above, and the second equalities are the associativity
$\alpha^x_{xxx}$ and  $(\alpha^x_{xxx})^{-1}$, respectively.

The same morphisms give a left rigidity structure when treated as
left birth and left death. Treat the objects $y$ and $\trivobject$
analogously by replacing $\phi$ with 1.

\begin{figure}

 \psset{unit=5mm}\begin{pspicture}[.3](0,0)(2,1)
\psline(0,0)(1,1) \psline(2,0)(1,1) \put (0.8,-1){$v^1_{xx}$}
\end{pspicture}\:\:\:,\:\:\:
\psset{unit=5mm}\begin{pspicture}[.3](0,0)(2,1) \psline(0,1)(1,0)
\psline(2,1)(1,0) \put(0.8,-1){$v^{x x}_1$}
\end{pspicture}\:\:\:,\:\:\:
\begin{pspicture}[.3](0,0)(2.5,3)
\psline(0,1)(1,2) \psline(2,1)(1,2) \psline(1,0)(0,1)
\psline(1,0)(2,1)
\end{pspicture} $=1$

\caption{Graphical notation of $v^1_{xx}$ and $v^{xx}_1$ and
property} \label{fig:basis}
\end{figure}

\begin{figure}

\psset{unit=4mm}
\begin{pspicture}[.3](0,0)(2,2)
\psbezier{<-}(0,0)(0,2)(2,2)(2,0)
\end{pspicture} $:=$
\begin{pspicture}[.3](0,0)(2,2)
\psline(0,0)(0,1) \psline(0,1)(1,2) \psline(1,2)(2,1)
\psline(2,0)(2,1)
\end{pspicture} $=:$ \begin{pspicture}[.3](0,0)(2,2)
\psbezier{->}(0,0)(0,2)(2,2)(2,0)
\end{pspicture} $\:\:\:,\:\:\:$
\begin{pspicture}[.3](0,0)(2,3)
\psbezier{<-}(0,2)(0,0)(2,0)(2,2)
\end{pspicture} $:= \frac{1}{\phi}$
\begin{pspicture}[.3](0,0)(2,2)
\psline(0,2)(0,1) \psline(0,1)(1,0) \psline(1,0)(2,1)
\psline(2,1)(2,2)
\end{pspicture} $=:$ \begin{pspicture}[.3](0,0)(2,3)
\psbezier{->}(0,2)(0,0)(2,0)(2,2)
\end{pspicture}

\begin{pspicture}[.3](0,0)(2,3)
\psbezier(0,1.5)(0,3)(2,3)(2,1.5)
\psbezier{->}(2,1.5)(2,0)(0,0)(0,1.5)
\end{pspicture} $= $
\begin{pspicture}[.3](0,0)(2,3)
\psbezier(0,1.5)(0,3)(2,3)(2,1.5)
\psbezier{<-}(2,1.5)(2,0)(0,0)(0,1.5)
\end{pspicture}
$= \frac{1}{\phi}$
\begin{pspicture}[.3](0,0)(2.5,3)
\psline(0,1)(1,2) \psline(2,1)(1,2) \psline(1,0)(0,1)
\psline(1,0)(2,1)
\end{pspicture}
$= \frac{1}{\phi}$.

\caption{Definitions of $b_x$ and $d_x$, and elementary properties}
\label{fig:b and d}
\end{figure}

\section{Proof of Theorem~\ref{thm_main} part \ref{main_nobraid}: the absence of braidings}\label{sec_nobraid}

The categories under consideration are known not to be braided (see
\cite{qa0503564}). However, once associativity matrices are known it
is in principle not difficult to classify braidings by direct
computation. In this section we perform this computation and show
that no braidings are possible.

A braiding consists of a natural family of isomorphisms $\{ c_{x,y}
: x \otimes y \rightarrow y \otimes x \} $ such that two hexagon
equalities hold:

$ (c_{x,y}\otimes z) \circ \alpha_{y,x,z}\circ (y \otimes c_{x,z})=
\alpha_{x,y,z}\circ c_{x,yz}\circ \alpha_{y,z,x} $ and

$ ((c_{y,x})^{-1}\otimes z) \circ \alpha_{y,x,z}\circ (y \otimes
(c_{z,x})^{-1}) = \alpha_{x,y,z}\circ
(c_{yz,x})^{-1}\circ\alpha_{y,z,x}
 $.

\begin{figure}
$ \xymatrix{&(yx)z \ar[rr]^{\alpha_{y,x,z}} && y(xz) \ar[rd]^{y
c_{x,z}} &\\(xy)z \ar[ru]^{c_{x,y}z} \ar[rd]^{\alpha_{x,y,z}}& & &
&y(zx)
 \\&x(yz) \ar[rr]^{c_{x,yz}} & &(yz)x \ar[ru]^{\alpha_{y,z,x}}&}$

$ \xymatrix{&(yx)z \ar[rr]^{\alpha_{y,x,z}} && y(xz) \ar[rd]^{y
c^{-1}_{z,x}} &\\(xy)z \ar[ru]^{c^{-1}_{y,x}z}
\ar[rd]^{\alpha_{x,y,z}}& & & &y(zx)
 \\&x(yz) \ar[rr]^{c^{-1}_{yz,x}} & &(yz)x \ar[ru]^{\alpha_{y,z,x}}&}$
\caption{Hexagon equalities} \label{fig:original hex}
\end{figure}

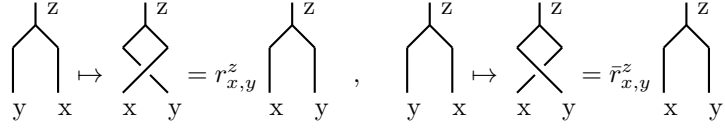
\begin{figure}\psset{unit=3mm}
$\begin{pspicture}[.3](0,0)(2,5) \psline(1,5)(1,4)
\psline(1,4)(0,3)\psline(1,4)(2,3) \psline(0,3)(0,1)
\psline(2,3)(2,1) \put(0,0){y}\put(2,0){x}\put(1.5,4.5){z}
\end{pspicture}\:\:
\mapsto \:\:
\begin{pspicture}[.3](0,0)(2,5) \psline(1,5)(1,4)
\psline(1,4)(0,3)\psline(1,4)(2,3) \psline(0,3)(2,1)
\psline[linecolor=white,linewidth=5pt](1.5,2.5)(0,1)
\psline(2,3)(0,1) \put(0,0){x}\put(2,0){y}\put(1.5,4.5){z}
\end{pspicture}\:\:
= r^z_{x,y}\:\:
\begin{pspicture}[.3](0,0)(2,5)
\psline(1,5)(1,4) \psline(1,4)(0,3)\psline(1,4)(2,3)
\psline(0,3)(0,1) \psline(2,3)(2,1)
\put(0,0){x}\put(2,0){y}\put(1.5,4.5){z}
\end{pspicture}$ $\;\;\;\;$, $\;\;\;\;$
$
\begin{pspicture}[.3](0,0)(2,5) \psline(1,5)(1,4)
\psline(1,4)(0,3)\psline(1,4)(2,3) \psline(0,3)(0,1)
\psline(2,3)(2,1) \put(0,0){y}\put(2,0){x}\put(1.5,4.5){z}
\end{pspicture}\:\:
\mapsto\:\:
\begin{pspicture}[.3](0,0)(2,5)
\psline(1,5)(1,4) \psline(1,4)(0,3)\psline(1,4)(2,3)
\psline(2,3)(0,1)
\psline[linecolor=white,linewidth=5pt](0.5,2.5)(2,1)
\psline(0,3)(2,1) \put(0,0){x}\put(2,0){y}\put(1.5,4.5){z}
\end{pspicture}\:\:
= \bar{r}^z_{x,y}\:\:
\begin{pspicture}[.3](0,0)(2,5)
\psline(1,5)(1,4) \psline(1,4)(0,3)\psline(1,4)(2,3)
\psline(0,3)(0,1) \psline(2,3)(2,1)
\put(0,0){x}\put(2,0){y}\put(1.5,4.5){z}
\end{pspicture}
 $
\caption{Isomorphisms $R^{z}_{x,y}$ and $\bar{R}^{z
}_{x,y}$}
\label{fig:braid}
\end{figure}

We define isomorphisms $R^z_{x,y}: V^z_{yx}\rightarrow V^z_{xy}$ by
$f \mapsto c_{x,y} \circ f $ and $\bar{R}^z_{x,y}:
V^z_{yx}\rightarrow V^z_{xy}$ by $f \mapsto (c_{x,y})^{-1} \circ f $
for any $f \in V^z_{yx}$. Figure \ref{fig:braid} shows the 1-dimensional case
where $r^z_{x,y}$ is nonzero and $\bar{r}^z_{x,y}=$
$(r^z_{y,x})^{-1}$. For higher dimensional spaces it can be
expressed as an invertible matrix, also denoted $r^z_{x,y}$ on the
canonically ordered basis as before.

With this linear isomorphism, the hexagon equations are equivalent
to the equations

$ \oplus_{s} R^{s}_{x,z} V^{t}_{ys} \circ \alpha^{t}_{y,x,z} \circ
\oplus_{s} R^{s}_{x,y} V^{t}_{sz} = \alpha^{t}_{y,z,x} \circ
\oplus_{s} V^{s}_{yz}R^{t}_{x,s} \circ \alpha^{t}_{x,y,z}$, and

\begin{figure}
\psset{unit=3mm,linewidth=0.5pt} $
\xymatrix{&\begin{pspicture}(0,0)(3,4) \psline(1,3.8)(1,3)
\psline(1,3)(0,2)\psline(1,3)(2,2)\psline(0.5,2.5)(1,2)
\psline(0,2)(1,1)\psline[linecolor=white,linewidth=5pt](0.7,1.7)(0.3,1.3)\psline(1,2)(0,1)
\psline(0,1)(0,0.6) \psline(1,1)(1,0.6)\psline(2,2)(2,0.6)
\put(-0.2,0){$x$}\put(0.8,0){$y$}\put(1.8,0){$z$}\put(0.1,2.8){$s$}\put(0.8,4){$t$}
\end{pspicture}
\ar[ld]_{\oplus_{s}R^{s}_{x,y}V^{t}_{sz} } &&
\begin{pspicture}(0,0)(3,4) \psline(1,3.8)(1,3)
\psline(1,3)(0,2)\psline(1,3)(2,2)\psline(1.5,2.5)(1,2)
\psline(0,2)(1,1)\psline[linecolor=white,linewidth=5pt](0.7,1.7)(0.3,1.3)\psline(1,2)(0,1)
\psline(0,1)(0,0.6) \psline(1,1)(1,0.6)\psline(2,2)(2,0.6)
\put(-0.2,0){$x$}\put(0.8,0){$y$}\put(1.8,0){$z$}\put(1.4,2.8){$s$}\put(0.8,4){$t$}
\end{pspicture}
\ar[ll]_{\alpha^{t}_{y,x,z}} &\\
\begin{pspicture}(0,0)(3,4) \psline(1,3.8)(1,3)
\psline(1,3)(0,2)\psline(1,3)(2,2) \psline(0,2)(0,0.6)
\psline(0.5,2.5)(1,2)\psline(1,2)(1,0.6)\psline(2,2)(2,0.6)
\put(-0.2,0){$x$}\put(0.8,0){$y$}\put(1.8,0){$z$}\put(0.1,2.8){$s$}\put(0.8,4){$t$}
\end{pspicture}
& &H^{t}_{x,y,z} & &
\begin{pspicture}(0,0)(3,4) \psline(1,3.8)(1,3)
\psline(1,3)(0,2)\psline(1,3)(2,2)\psline(1.5,2.5)(1,2)
\psline(0,2)(1,1)\psline(1,2)(2,1)\psline[linecolor=white,linewidth=5pt](1.6,1.8)(0.2,1.1)
\psline(2,2)(0,1) \psline(0,1)(0,0.6)
\psline(1,1)(1,0.6)\psline(2,1)(2,0.6)
\put(-0.2,0){$x$}\put(0.8,0){$y$}\put(1.8,0){$z$}\put(1.4,2.8){$s$}\put(0.8,4){$t$}
\end{pspicture}
\ar[lu]_{\oplus_{s} R^{s}_{x,z}V^{t}_{ys}}
\ar[ld]^{\alpha^{t}_{y,z,x}}
 \\&
 \begin{pspicture}(0,0)(3,4)
\psline(1,3.8)(1,3) \psline(1,3)(0,2)\psline(1,3)(2,2)
\psline(0,2)(0,0.6)
\psline(1.5,2.5)(1,2)\psline(1,2)(1,0.6)\psline(2,2)(2,0.6)
\put(-0.2,0){$x$}\put(0.8,0){$y$}\put(1.8,0){$z$}\put(1.4,2.8){$s$}\put(0.8,4){$t$}
\end{pspicture}
\ar[lu]^{\alpha^{t}_{x,y,z}} & &
 \begin{pspicture}(0,0)(3,4) \psline(1,3.8)(1,3)
\psline(1,3)(0,2)\psline(1,3)(2,2)\psline(0.5,2.5)(1,2)
\psline(0,2)(1,1)\psline(1,2)(2,1)\psline[linecolor=white,linewidth=5pt](1.6,1.8)(0.2,1.1)
\psline(2,2)(0,1) \psline(0,1)(0,0.6)
\psline(1,1)(1,0.6)\psline(2,1)(2,0.6)
\put(-0.2,0){$x$}\put(0.8,0){$y$}\put(1.8,0){$z$}\put(0.1,2.8){$s$}\put(0.8,4){$t$}
\end{pspicture}
\ar[ll]_{\oplus_{s} V^{s}_{yz}R^{t}_{x,s}}&}$

$ \xymatrix{& \begin{pspicture}(0,0)(3,5) \psline(1,3.8)(1,3)
\psline(1,3)(0,2)\psline(1,3)(2,2)\psline(0.5,2.5)(1,2)
\psline(1,2)(0,1)\psline[linecolor=white,linewidth=5pt](0.2,1.8)(0.7,1.2)\psline(0,2)(1,1)
\psline(0,1)(0,0.6) \psline(1,1)(1,0.6)\psline(2,2)(2,0.6)
\put(-0.2,0){$x$}\put(0.8,0){$y$}\put(1.8,0){$z$}\put(0.1,2.8){$s$}\put(0.8,4){$t$}
\end{pspicture}
\ar[ld]_{\oplus_{s}\bar{R}^{s}_{x,y}V^{t}_{sz} } &&
\begin{pspicture}(0,0)(3,5) \psline(1,3.8)(1,3)
\psline(1,3)(0,2)\psline(1,3)(2,2)\psline(1.5,2.5)(1,2)
\psline(1,2)(0,1)\psline[linecolor=white,linewidth=5pt](0.2,1.8)(0.7,1.2)\psline(0,2)(1,1)
\psline(0,1)(0,0.6) \psline(1,1)(1,0.6)\psline(2,2)(2,0.6)
\put(-0.2,0){$x$}\put(0.8,0){$y$}\put(1.8,0){$z$}\put(1.4,2.8){$s$}\put(0.8,4){$t$}
\end{pspicture}\ar[ll]_{\alpha^{t}_{y,x,z}} &\\
\begin{pspicture}(0,0)(3,4)
\psline(1,3.8)(1,3) \psline(1,3)(0,2)\psline(1,3)(2,2)
\psline(0,2)(0,0.6)
\psline(0.5,2.5)(1,2)\psline(1,2)(1,0.6)\psline(2,2)(2,0.6)
\put(-0.2,0){$x$}\put(0.8,0){$y$}\put(1.8,0){$z$}\put(0.1,2.8){$s$}\put(0.8,4){$t$}
\end{pspicture} & & \bar{H}^{t}_{x,y,z}& &
\begin{pspicture}(0,0)(3,4)\psline(1,3.8)(1,3)
\psline(1,3)(0,2)\psline(1,3)(2,2)\psline(1.5,2.5)(1,2)
\psline(2,2)(0,1)\psline[linecolor=white,linewidth=5pt](0.2,1.8)(0.7,1.2)
\psline[linecolor=white,linewidth=5pt](1.2,1.8)(1.7,1.2)
\psline(0,2)(1,1)\psline(1,2)(2,1) \psline(0,1)(0,0.6)
\psline(1,1)(1,0.6)\psline(2,1)(2,0.6)
\put(-0.2,0){$x$}\put(0.8,0){$y$}\put(1.8,0){$z$}\put(1.4,2.8){$s$}\put(0.8,4){$t$}
\end{pspicture}\ar[lu]_{\oplus_{s}
\bar{R}^{s}_{x,z}V^{t}_{ys}}\ar[ld]^{\alpha^{t}_{y,z,x}}
 \\&\begin{pspicture}(0,0)(3,4)
\psline(1,3.8)(1,3) \psline(1,3)(0,2)\psline(1,3)(2,2)
\psline(0,2)(0,0.6)
\psline(1.5,2.5)(1,2)\psline(1,2)(1,0.6)\psline(2,2)(2,0.6)
\put(-0.2,0){$x$}\put(0.8,0){$y$}\put(1.8,0){$z$}\put(1.4,2.8){$s$}\put(0.8,4){$t$}
\end{pspicture}
\ar[lu]^{\alpha^{t}_{x,y,z}}
 \ar[lu]^{\alpha^{t}_{x,y,z}} & &
 \begin{pspicture}(0,0)(3,4) \psline(1,3.8)(1,3)
\psline(1,3)(0,2)\psline(1,3)(2,2)\psline(0.5,2.5)(1,2)
\psline(2,2)(0,1)\psline[linecolor=white,linewidth=5pt](0.2,1.8)(0.7,1.2)
\psline[linecolor=white,linewidth=5pt](1.2,1.8)(1.7,1.2)
\psline(0,2)(1,1)\psline(1,2)(2,1) \psline(0,1)(0,0.6)
\psline(1,1)(1,0.6)\psline(2,1)(2,0.6)
\put(-0.2,0){$x$}\put(0.8,0){$y$}\put(1.8,0){$z$}\put(0.1,2.8){$s$}\put(0.8,4){$t$}
\end{pspicture} \ar[ll]_{\oplus_{s} V^{s}_{yz}\bar{R}^{t}_{x,s}}&}$
\caption{Equivalent hexagon equalities} \label{fig:hex}
\end{figure}
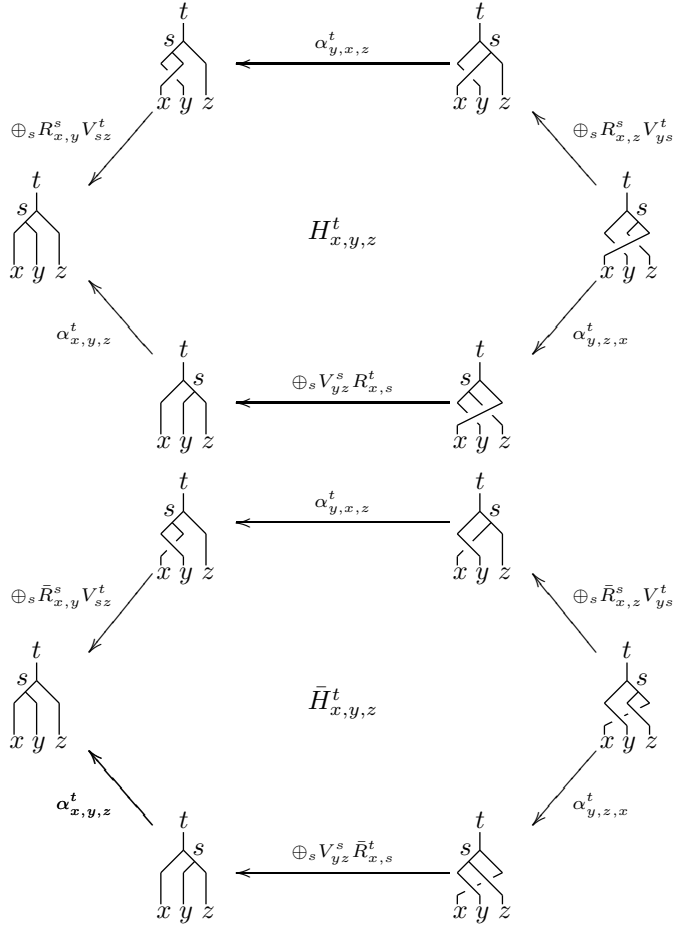

$ \oplus_{s} \bar{R}^{s}_{x,z} V^{t}_{ys} \circ \alpha^{t}_{y,x,z}
\circ \oplus_{s} \bar{R}^{s}_{x,y} V^{t}_{sz} = \alpha^{t}_{y,z,x}
\circ \oplus_{s} V^{s}_{yz}\bar{R}^{t}_{x,s} \circ
\alpha^{t}_{x,y,z}$, which we still call hexagon equations, referred
to as $H^{t}_{x,y,z}$ and $\bar{H}^{t}_{x,y,z}$, respectively. These are illustrated graphically in
Figure~\ref{fig:hex}).

We show the absence of a braiding by assuming the existence and
deriving a contradiction.

We need five 2-dimensional hexagon equations as follows:

$H^{x}_{y,x,x} : R^{x}_{y,x}\otimes I_2 \circ \alpha^{x}_{x,y,x}
\circ R^{x}_{y,x}\otimes I_2 = \alpha^{x}_{x,x,y} \circ I_2\otimes
R^{x}_{y,x} \circ \alpha^{x}_{y,x,x}  $

$\bar{H}^{x}_{y,x,x} : \bar{R}^{x}_{y,x}\otimes I_2 \circ
\alpha^{x}_{x,y,x} \circ \bar{R}^{x}_{y,x} \otimes I_2 =
\alpha^{x}_{x,x,y} \circ  I_2\otimes \bar{R}^{x}_{y,x} \circ
\alpha^{x}_{y,x,x}  $

$H^{x}_{x,y,x} : R^{x}_{x,x}\otimes 1 \circ \alpha^{x}_{y,x,x}
\circ R^{x}_{x,y} \otimes I_2 = \alpha^{x}_{y,x,x} \circ 1\otimes
R^{x}_{x,x} \circ \alpha^{x}_{x,y,x}  $

$\bar{H}^{x}_{x,y,x} : \bar{R}^{x}_{x,x}\otimes 1 \circ
\alpha^{x}_{y,x,x} \circ \bar{R}^{x}_{x,y} \otimes I_2 =
\alpha^{x}_{y,x,x} \circ  1\otimes \bar{R}^{x}_{x,x} \circ
\alpha^{x}_{x,y,x}  $

$H^{1}_{x,x,x} : R^{x}_{x,x}\otimes 1 \circ \alpha^{1}_{x,x,x}
\circ R^{x}_{x,x} \otimes 1 = \alpha^{1}_{x,x,x} \circ  I_2\otimes
R^{1}_{x,x} \circ \alpha^{1}_{x,x,x}  $

These are of the following forms, respectively:

$ (r^{x}_{y,x})^{2}\left[\begin{array}{cc}1
&0\\0&-1\end{array}\right] =r^{x}_{y,x}\left[\begin{array}{cc}0
&b\\1/b & 0\end{array}\right]\left[\begin{array}{cc}0 &1\\1&
0\end{array}\right]$

$ (r^{x}_{x,y})^{-2}\left[\begin{array}{cc}1
&0\\0&-1\end{array}\right]
=(r^{x}_{x,y})^{-1}\left[\begin{array}{cc}0 &b\\1/b &
0\end{array}\right]\left[\begin{array}{cc}0 &1\\1&
0\end{array}\right]$

$r^{x}_{x,y}\left[\begin{array}{cc}k
&l\\m&n\end{array}\right]\left[\begin{array}{cc}0 &1\\1&
0\end{array}\right] = \left[\begin{array}{cc}0 &1\\1&
0\end{array}\right] \left[\begin{array}{cc}k
&l\\m&n\end{array}\right] \left[\begin{array}{cc}1
&0\\0&-1\end{array}\right] $

$ (r^{x}_{y,x})^{-1} \left[\begin{array}{cc}k &l\\m&n
\end{array}\right]^{-1} \left[\begin{array}{cc}0 &1\\1&
0\end{array}\right] = \left[\begin{array}{cc}0 &1\\1&
0\end{array}\right] \left[\begin{array}{cc}k
&l\\m&n\end{array}\right]^{-1} \left[\begin{array}{cc}1
&0\\0&-1\end{array}\right]   $

$d \left[\begin{array}{cc}k
&l\\m&n\end{array}\right]\left[\begin{array}{cc}1 &b\\1&
-b\end{array}\right]\left[\begin{array}{cc}k
&l\\m&n\end{array}\right] = d^2 r^1_{x,x}\left[\begin{array}{cc}1
&b\\1& -b\end{array}\right]^2 $

where $\left[\begin{array}{cc}k &l\\m&n
\end{array}\right]$ represents the matrix $r^{x}_{x,x}$.

From the first four equations, we get $r^{x}_{y,x}=b$, $r^{x}_{x,y}=
1/b$, $-n=r^{x}_{x,y} k$, $m=r^{x}_{x,y} l$,
 $-n=r^{x}_{y,x} k$, which imply $k=n=0$ since $r^x_{x,y}
\neq r^x_{y,x}$ as above. Now from the final one we get $l^2 = d
r^1_{x,x}(b+1)$ and $-blm=d r^1_{x,x}(1+b)$, and the later equality
means $l^2 =-dr^1_{x,x}(1+b)$ by substituting $m=r^{x}_{x,y} l$. We
get easily a contradiction for either case $b= \pm i$.

\section{Pivotal structures and sphericity}\label{sec_pivotalprelim}
Let $\cat{C}$ be a rigid monoidal category. A {\em pivotal structure} for $\cat{C}$ is a monoidal natural isomorphism $\pi$ from $**$ to $Id$. A {\em strict pivotal structure} is a pivotal structure which is the identity. In a pivotal monoidal category, the {\em right trace} $tr_r$ of an endomorphism $f:x \to x$ is given by $tr_r(f)=b_x \circ (f \otimes Id_{x^*}) \circ (\pi_x^{-1} \otimes Id_{x^*}) \circ d_{x^*} \in End(\trivobject) \cong \fld{C}$. The {\em left trace} $tr_l$ is given by $tr_l(f)=b_{x^*} \circ (f^* \otimes Id_{x^{**}}) \circ ((\pi_x)^* \otimes Id_{x^{**}}) \circ d_{x^{**}}$. A pivotal monoidal category is {\em spherical} if $tr_r=tr_l$.

Pivotal structures may not be unique. For example, in a fusion category with object types given by a finite group $G$, group multiplication as tensor product and trivial associativity matrices, any group homomorphism $G \to \fld{C}$ induces a pivotal structure. Furthermore, pivotal structures depend on choices of rigidity. However, if one chooses a new rigidity structure with $b'_x = c b_x$ and $d'_x = c^{-1} d_x$, then $\pi'_x = c^{-1} \pi_x$ gives a new pivotal structure $\pi'$ inducing the same traces as $\pi$.

For a strictified skeletal fusion category, we shall assume the rigidity structures described in Section~\ref{sec_strictskel}. Then $**$ is an object fixing monoidal endofunctor. The isomorphisms $J_{a,b}:a^{**} \otimes b^{**} \to (a \otimes b)^{**}$ associated with $**$ considered as a monoidal functor may be taken to be the identity on $a \otimes b$. If such a category has a pivotal structure $\pi$, it must take the following form: for each strand $(x)$, there is a scalar $t_x$ such that $\pi_x = t_x Id_x$, and $\pi_{(x_1,\ldots,x_n)} = t_{x_1} \ldots t_{x_n} Id_{(x_1, \ldots x_n)}$. Then for all sequences $(x_1, \ldots, x_m)$ and $(y_1, \ldots, y_n)$, and all $f:(x_1, \ldots, x_m) \to (y_1, \ldots, y_n)$, $f^{**}=t_{x_1} \ldots t_{x_m} t_{y_1}^{-1} \ldots t_{y_n}^{-1} f$. This implies that, in particular, $t_\trivobject=1$.

Writing out the diagrams for $b_x^{**}$ and $d_x^{**}$ and applying rigidity gives that $b_x^{**}=b_x$ and $d_x^{**}=d_x$. Thus one must have $t_x t_{x^*} = 1$, and for self-dual strands, $t_x = \pm 1$; in this case $t_x$ is called the {\em Frobenius-Schur indicator} for $x$.

Furthermore, in a strictified skeleton the left and right trace on a strand may be rewritten as follows:
\[tr_r(f)=t_x^{-1} b_x \circ (f \otimes Id_{x^*}) \circ d_{x^*},\]
\[tr_l(f)=t_x b_{x^*} \circ (Id_{x^*} \otimes f) \circ d_x.\]

\begin{lmm}\label{lmm_selfdual}
Every pivotal fusion category with self-dual simple objects is spherical.
\end{lmm}
\begin{proof}
The result holds since if $x$ is a self dual strand then $b_x=b_{x^*}$, $d_x=d_{x^*}$ and $t_x=t_x^{-1}$. Thus for each $f:x \to x$ we have $tr_r(f)=tr_l(f)$.
\end{proof}

Kitaev has shown in \cite{kitaev} that every unitary category admits a spherical structure. A more general property called {\em pseudo-unitarity} is shown in \cite{etingof_nikshych_ostrik} to guarantee a spherical structure. However, it is not known if every fusion category admits a pivotal or spherical structure.

For arbitrary fusion categories, one has that $**** \cong Id$. This was shown in \cite{etingof_nikshych_ostrik}, using an analog of Radford's formula for $S^4$ for representation categories of weak Hopf algebras, which was developed in \cite{nikshych}. The following theorem shows that, in a strictified skeletal fusion category, a convenient choice of rigidity makes $****$ the identity on the nose. Extending the result to general fusion categories via natural equivalences gives an elementary proof that $**** \cong Id$.

\begin{thrm}
In a strictified skeletal fusion category, there is a choice of rigidity structures such that $****=Id$.
\end{thrm}
\begin{proof}
The functor $****$ is the identity on $(2)$-stranded morphisms by rigidity; it suffices to prove the result for $(2,1)$ stranded morphisms. Let $V=V_{x y}^z$ be a $(2,1)$ stranded morphism space with a basis $\{v_i\}$, and let $\{w_i\}$ be an algebraically dual basis for the space $W=V_z^{x y}$, in the sense that $w_i \circ v_j=\delta_{ij} Id_z$. For any simple object $z$, define the {\em right pseudo-trace} $ptr_r$ of an endomorphism $f:z \to z$ by $ptr_r(f)=b_{z} \circ (f \otimes Id_{z^*}) \circ d_{z^*}$, and the {\em left pseudo-trace} $ptr_l$ by $ptr_l(f) = b_{z^*} \circ (Id_{z^*} \otimes f) \circ d_{z}$. This definition is possible because $**$ is the identity on objects. Scale rigidity morphisms if necessary so that for any strand $z$, $ptr_r(Id_z)=ptr_l(Id_z)$. Because $d_z$ and $b_{z^*}$ are nonzero elements of one dimensional algebraically dual morphism spaces, $ptr_r(Id_z) \ne 0$. One may now exchange left pseudo-traces for right pseudo-traces, just like with traces in a graphical calculus for a spherical category. 

Figure~\ref{fig_quadrupledual} gives the proof. On the left side, bending arms and pseudo-sphericity implies that the algebraic dual basis of the basis $\{w_i^{**}\}$ is $\{\leftexp{**}{v_i}\}$. However, on the right side the functoriality of the double dual implies that the algebraic dual basis of $\{w_i^{**}\}$ is $\{v_i^{**}\}$. Since the left and right double dual are inverse functors, $****$ is the identity.
\end{proof} 

\standardfigsized{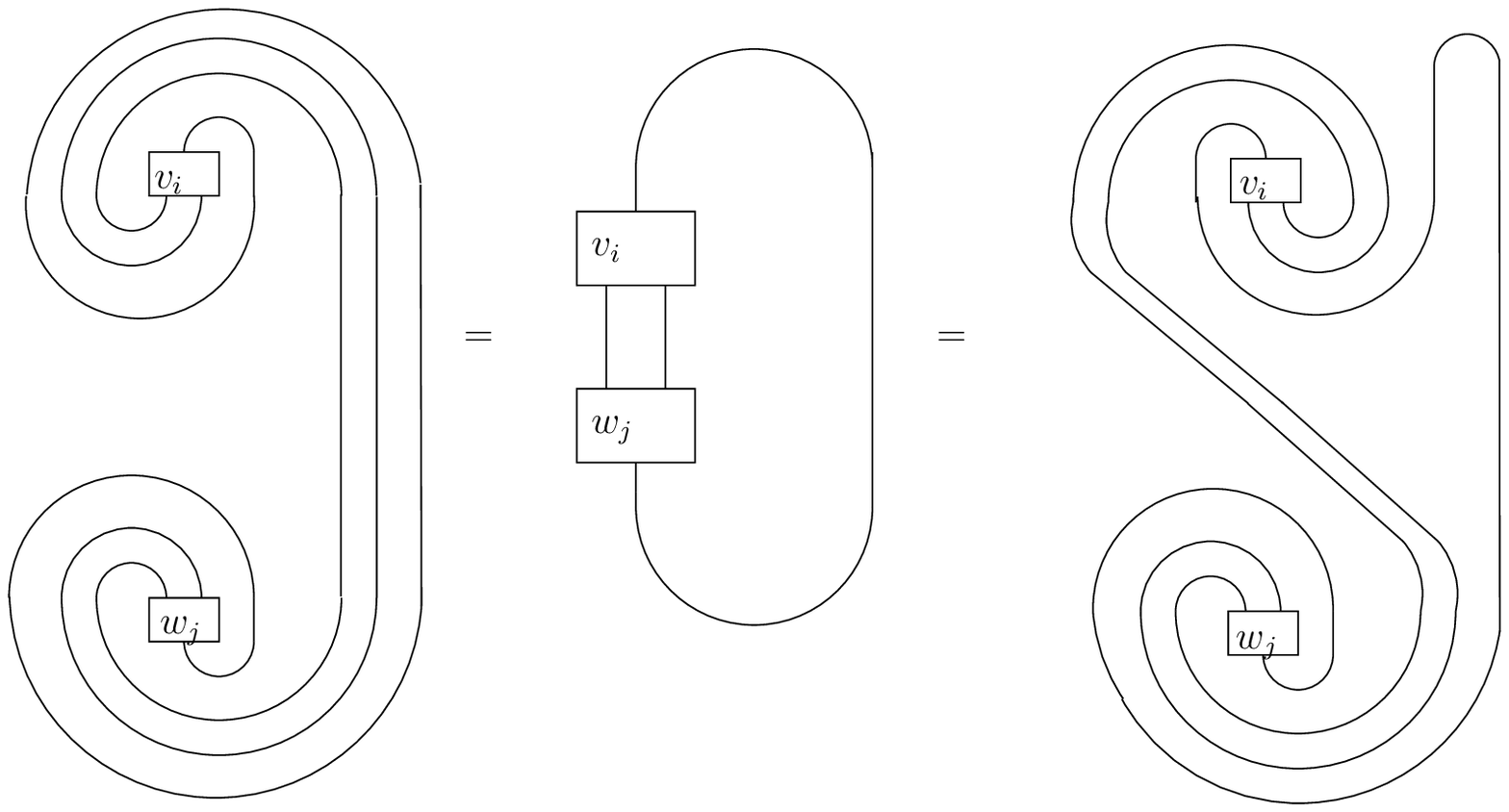}{In a strictified skeletal fusion category, with the right choice of rigidity structures the quadruple dual is the identity.}{scale=.7}

Even if a category admits a pivotal structure it is not known whether it admits a spherical pivotal structure. Pictorial considerations do not readily provide an answer. It is possible, however, to partially describe what a pivotal strictified skeleton which did not admit a spherical structure would look like.

Let $\cat{C}$ be a pivotal strictified skeletal fusion category which does not admit a spherical structure. Choose rigidity morphisms which give a pseudo-spherical structure as above, and a matching pivotal structure. Then for any object $x$, one has the following:
\[\frac{tr_l(Id_x)}{tr_r(Id_x)}= \frac{t_x ptr_l(x)}{t_x^{-1} ptr_r(x)}=t_x^2.\]
Therefore, $\cat{C}$ is spherical iff there exists a pivotal structure such that all of the $t_x$ are $\pm 1$. Thus there must be some strand $x$ such that $t_x \neq \pm 1$.

For strands $u$ and $v$, $u \otimes v$ has a nontrivial morphism to some object $w$, and $t_u t_v (t_w)^{-1}=\pm 1$, since $****=Id$. Thus the set of scalars $t$ and their additive inverses forms a finite subgroup $G$ of $\fld{C}$. Note that we can apply any group homomorphism that preserves $\pm 1$ to the set of scalars $t$ and get a new set of scalars $t'$ which also gives a pivotal structure. At least one product $t_u t_v (t_w)^{-1}$ must be equal to $-1$, or else we could apply the trivial homomorphism to the set of scalars $t$ to get a new pivotal structure with $t'_u=1$ for all strands $u$, which would make $\cat{C}$ spherical.

Every finite subgroup of $\fld{C}$ is a cyclic group of roots of unity. We have $|G|=2k$ for some $k$, and since $\cat{C}$ is not spherical, $|G| \ge 4$. Using a homomorphism which preserves $-1$ we may switch to a new pivotal structure which gives $|G|=2^k$ for some $k$, where $k \ge 2$ to contradict sphericity. Pick an object $v$ with $t_v^2=-1$. Then $v$ is not self dual, and for a simple summand $w$ in $v \otimes v$, one has $w \ne \trivobject$ and $t_w^2=1$. Therefore, $\cat{C}$ has at least four objects, $v$, $v^*$, $w$ and $\trivobject$. The set of objects $u$ such that $t_u^2=1$ generates a spherical subcategory $\cat{C}'$ with at least two simple objects, and missing at least two.

\begin{lmm}
Any fusion category which is pivotal but admits no spherical structure contains at least five simple object types.
\end{lmm}
\begin{proof}
Assume that $\cat{C}$ has four simple object types, $\trivobject$, $w$, $v$ and $v'$ as above. Then $\cat{C}'$ has two simple object types, and by the classification of fusion categories with two simple object types in \cite{ostrik2}, its fusion rules are given by $w \otimes w \cong n w \oplus \trivobject$, where $n \in \{0,1\}$.  The pivotal structure places limitations on the fusion rules, for example $v \otimes w \cong a v \oplus b v'$ for some $a$ and $b$ in $\fld{N}$. An easy calculation shows that $\cat{C}$ admits only one associative fusion ring, in which objects and tensor products are given by the group $\fld{Z}_4$. Any such category is pseudo-unitary and therefore spherical, as described in \cite{etingof_nikshych_ostrik}, which contradicts the assumption. Therefore, a pivotal fusion category which can't be made spherical must have at least five simple object types.
\end{proof}
\section{Proof of Theorem~\ref{thm_main} part \ref{main_spherical}: spherical structure calculations}\label{sec_pivotal}
In this section we explicitly compute pivotal structures for the categories found in Section~\ref{sec_tensor}. Since these categories have self dual simple objects, Lemma~\ref{lmm_selfdual} implies that they are spherical.

It is not hard to determine whether or not a fusion category is pivotal once a set of associativity matrices is known. One way is to perform the calculations directly using the associativity matrices, but there is an easier calculation. In order to explain this calculation, it is convenient to extend the definition of composition of morphisms over extra-categorical direct sums of morphism spaces. Suppose $f \in Mor(a,b)$ and $g \in Mor(c,d)$. Define $f \circ g$ as usual if $b=c$, and $f \circ g = 0 \in Mor(a,d)$ otherwise. Extend this definition over direct sums of morphism spaces, distributing composition over direct sum.

Given a strictified skeletal fusion category $C$ and a set of associativity matrices, choose bases for the $(2,0)$ and $(2,1)$-stranded morphism spaces compatible with the associativity matrices and choose rigidity so that for each strand $x$, the basis element for $V_{x^* x}^1$ is $d_x$. Define morphisms $b=\oplus_{x} b_x$, $d=\oplus_x d_x$, and $I=\oplus_x Id_x$, taking sums over the strands.

Then $B$ acts on $\bigoplus_{x,y,z} V_{x y}^z$ as follows:
\[B(f)=(I \otimes I \otimes b) \circ (I \otimes f \otimes I) \circ (d \otimes I)\]

For a single $(2,1)$-stranded morphism space, this action amounts to ``bending arms''. The cube of $B$ is the double dual. The action of $B$ on a morphism $f \in V_{x y}^z$ is given by the associativity matrix $a_{z^*,x,y}^\trivobject$, since $(Id_{z^*} \otimes f) \circ d_{z} = (g \otimes id_{y}) \circ d_{y}$ for some $g \in V_{z^* x}^{y^*}$ implies that $B(f)=(Id_{z^*} \otimes Id_x \otimes b_y) \circ (Id_{z^*} \otimes f \otimes Id_{y^*}) \circ (d_{z} \otimes Id_{y^*}) = (Id_{z^*} \otimes Id_x \otimes b_{y}) \otimes (g \otimes Id_{y} \otimes Id_{y^*}) \circ (d_{y} \otimes Id_{y^*}) = g$ by rigidity. For the fusion rules at hand, the matrix for $B$ is as follows:

$$ \begin{array}{l|ccccc}
          & v_1 & v_2 & v_{x x}^y & v_{y x}^x & v_{x y}^x \\
\hline
v_1       & (a_{x x x}^\trivobject)_{1,1} & (a_{x x x}^\trivobject)_{1,2} & 0 & 0 & 0 \\
v_2       & (a_{x x x}^\trivobject)_{2,1} & (a_{x x x}^\trivobject)_{2,2} & 0 & 0 & 0 \\
v_{x x}^y & 0 & 0 & 0 & a_{y x x}^\trivobject & 0 \\
v_{y x}^x & 0 & 0 & 0 & 0 & a_{x y x}^\trivobject \\
v_{x y}^x & 0 & 0 & a_{x x y}^\trivobject & 0 & 0
\end{array}
$$
For all of the solutions given in Section~\ref{sec_tensor}, $B^3$ is the identity matrix, so the corresponding strictified categories have a strict pivotal structure. Non-strict pivotality would mean that $B^3$ is a diagonal matrix with eigenvalues determined by a family of invertible scalars $t$, coherent as described in Section~\ref{sec_pivotalprelim}.

\appendix

\section{Associativity matrices}\label{sec_appendix}

In this section, we give explicit associativity matrices for the categorical realization given in Section~\ref{sec_tensor}.

$a^{y}_{y,y,y}=$ $a^{x}_{x,y,y}=$ $a^{x}_{y,y,x}=$
$a^{1}_{x,y,x}=$ $a^{1}_{x,x,y}=$ $a^{y}_{x,x,y}=$
$a^{1}_{y,x,x}=$ $a^{y}_{y,x,x}=1$,

$a^{y}_{x,y,x}=$ $a^{x}_{y,x,y}= -1$,

$a^{x}_{x,y,x}= \left[\begin{array}{cc}1&0\\0&-1\end{array}
\right]$

$a^{x}_{x,x,y}= \left[\begin{array}{cc}0&i\\-i&0\end{array}
\right]$

$a^{x}_{y,x,x}= \left[\begin{array}{cc}0&1\\1&0\end{array}
\right]$

$a^{1}_{x,x,x}= \frac{1}{\sqrt{2}} e^{7 \pi i
 /12} \left[\begin{array}{cc}1&i\\1&-i\end{array}
\right]$

$a^{y}_{x,x,x}= \frac{1}{\sqrt{2}} e^{7 \pi i
 /12} \left[\begin{array}{cc}i&1\\-i&1\end{array}
\right]$

$a^{x}_{x,x,x}=$ $\left[ \begin{smallmatrix}
\frac{-1+\sqrt{3}}{2}&\frac{-1+\sqrt{3}}{2}& \frac{1- \sqrt{3}}{4}
e^{\pi i /6}&\frac{1- \sqrt{3}}{4} e^{2\pi i /3}
&\frac{1- \sqrt{3}}{4} e^{2\pi i /3}&\frac{1- \sqrt{3}}{4} e^{\pi i /6} \\
\frac{-1+\sqrt{3}}{2}&\frac{1-\sqrt{3}}{2}& \frac{1- \sqrt{3}}{4}
e^{\pi i /6}&\frac{1- \sqrt{3}}{4} e^{2\pi i /3}&-\frac{1-
\sqrt{3}}{4} e^{2\pi i /3}&-\frac{1- \sqrt{3}}{4} e^{\pi i /6}
\\1&1&-\frac{1}{2} (e^{\pi i /6} -1)&\frac{1}{2} e^{5\pi i /6}&\frac{1}{2} (e^{-\pi i /3} + i)&\frac{1}{2} e^{\pi i /3}
\\1&1&\frac{1}{2} e^{\pi i /3}&\frac{1}{2} (e^{-\pi i /3} + i)&\frac{1}{2} e^{5\pi i /6}&-\frac{1}{2} (e^{\pi i /6} -1)
\\1&-1&-\frac{1}{2} (e^{\pi i /6} -1)&\frac{1}{2} e^{5\pi i /6}&-\frac{1}{2} (e^{-\pi i /3} +
i)&-\frac{1}{2} e^{\pi i /3}
\\-1&1&-\frac{1}{2} e^{\pi i /3}&-\frac{1}{2} (e^{-\pi i /3} + i)&\frac{1}{2} e^{5\pi i /6}&-\frac{1}{2} (e^{\pi i /6} -1)
\end{smallmatrix}
\right]$

\bibliographystyle{chicago}
\bibliography{references}

\begin{thebibliography}{10}

\bibitem{buchberger}
Bruno Buchberger.
\newblock A theoretical basis for the reduction of polynomials to canonical
  forms.
\newblock {\em SIGSAM Bull.}, 10(3):19--29, 1976.

\bibitem{etingof_nikshych_ostrik}
Pavel Etingof, Dmitri Nikshych, and Viktor Ostrik.
\newblock On fusion categories.
\newblock {\em Ann. of Math.}, 162(2):581--642, 2005.

\bibitem{kerler}
J.~Frohlich and T.~Kerler.
\newblock {\em Quantum groups, quantum categories, and quantum field theory},
  chapter~4.
\newblock Number 1542 in Lecture Notes in Mathematics. 1993.

\bibitem{kassel}
Christian Kassel.
\newblock {\em Quantum Groups}.
\newblock Springer-Verlag, 1995.

\bibitem{kazhdan_wenzl}
Kazhdan and Hans Wenzl.
\newblock Reconstructing monoidal categories.
\newblock {\em Adv. Soviet Math.}, 16:111--136, 1993.

\bibitem{kitaev}
Alexei Kitaev.
\newblock Anyons in an exactly solved model and beyond.
\newblock {\em Annals of Physics}, 321(1):2--111, 2006.

\bibitem{maclane}
Saunders Mac~Lane.
\newblock {\em Categories for the Working Mathematician, Second Edition}.
\newblock Springer-Verlag, 1978.

\bibitem{nikshych}
Dmitri Nikshych.
\newblock On the structure of weak hopf algebras.
\newblock {\em And. Math.}, 170:257--286, 2002.

\bibitem{qa0503564}
Victor Ostrik.
\newblock Pre-modular categories of rank 3.
\newblock math.CT/0507349.

\bibitem{ostrik2}
Victor Ostrik.
\newblock Fusion categories of rank 2.
\newblock {\em Math. Res. Lett.}, 10(2-3):177--183, 2003.

\bibitem{rowell_stong_wang}
Eric Rowell, Richard Stong, and Zhenghan Wang.
\newblock in preparation.

\bibitem{tambara_yamagami}
Daisuke Tambara and Shigeru Yamagami.
\newblock Tensor categories with fusion rules of self-duality for finite
  abelian groups.
\newblock {\em J. Algebra}, 209:692--707, 1998.

\bibitem{wenzl_tuba}
Hans Wenzyl and Imre Tuba.
\newblock On braided tensor categories of type bcd.
\newblock {\em J. Reine. Angew. Math.}, 581:31--69, 2005.

\end{thebibliography}

\end{document}